\documentclass[onefignum,onetabnum]{siamart220329}
\usepackage{preambleCommands}
\setuptodonotes{disable}	 
\newcommand{\loli}{\emph{locality limitation}\xspace}
\headers{Locality Limitation of CG}{U. R{\"u}de}
\theoremstyle{remark}
\newsiamremark{remark}{Remark}
\crefname{remark}{Remark}{Remarks}
\Crefname{remark}{Remark}{Remarks}
\begin{document}
\title{Conjugate Gradient Methods are Not Efficient  \\[1ex]
	Experimental Study of the Locality Limitation
	}

	\author{U. R\"ude\thanks{
	(\email{ulrich.ruede@fau.de}) \;
    	Department of Applied Mathematics, VSB-Technical University of Ostrava, Czech Republic;\;
	CERFACS, Toulouse, France; \;
	Lehrstuhl für Systemsimulation, Friedrich-Alexander Universt{\"a}t Erlangen-Nürnberg.} }
%
%
\maketitle
\begin{abstract}
The convergence of the Conjugate Gradient method is subject to a \loli, 
which imposes a lower bound on the number of iterations required before a 
qualitatively accurate approximation can be obtained. 
This limitation originates from the restricted transport of information in the graph induced 
by the sparsity pattern of the system matrix. 
In each iteration, information from the right-hand side can propagate only across directly connected graph nodes. 
The diameter of this graph therefore determines a minimum number of
iterations that is necessary to achieve an acceptable level of accuracy.
\end{abstract}
%
%


\section{Introduction} 
This report presents results from an experimental study investigating the influence of a \loli 
on the convergence behavior of the \ac{CG} method and related iterative solvers. 
The \loli phenomenon arises from the inherently local nature of data transport in iterative algorithms. 
In particular, information propagation is governed by the graph associated with the system matrix: 
during each matrix–vector multiplication, data is exchanged only between neighboring nodes connected by edges of this graph.
As a consequence, each iteration of the CG method can propagate information relevant to the construction of the global solution by at most one graph edge. 

This implies that, independent of spectral considerations, 
a minimal number of iterations may be required before a globally accurate solution can be attained. 
The present study analyzes this effect using sparse matrices derived from the discretization of 
one- and two-dimensional boundary value problems. 
The test cases are deliberately constructed to isolate and illustrate the impact of the \loli 
on convergence.
In addition, several variations of the basic setup are examined in order to assess the generality and the significance
of the observed phenomena.%
\begin{algorithm}[ht]
\caption{\ac{CG}}
\begin{algorithmic}[1]
\REQUIRE matrix $\Am$ \ac{SPD}, vector $\bv$, initial guess $\uv_0$, tolerance $\tol$
\STATE $\rv_0 \leftarrow \bv - \Am \uv_0$ \label{CG-residual0}
\STATE $\pv_0 \leftarrow \rv_0$
\STATE $k \leftarrow 0$
\WHILE{$\| \rv_k\| > \tol$}
    \STATE $\alpha_k \leftarrow \dfrac{\rv_k^T \rv_k}{\pv_k^T \Am \pv_k}$
    \STATE $\uv_{k+1} \leftarrow \uv_k + \alpha_k \pv_k$
    \STATE $\rv_{k+1} \leftarrow \rv_k - \alpha_k \Am \pv_k$ \label{CG-residual}
    \IF{$\| \rv_{k+1}\| \leq \tol$} \label{exit-line}
        \STATE \textbf{break}
    \ENDIF
    \STATE $\beta_k \leftarrow \dfrac{ \rv_{k+1}^T \rv_{k+1}}{\rv_k^T \rv_k}$
    \STATE $\pv_{k+1} \leftarrow \rv_{k+1} + \beta_k \pv_k$
    \STATE $k \leftarrow k + 1$
\ENDWHILE
\RETURN $\uv_{k+1}$
\end{algorithmic}
\label{alg:cg}
\end{algorithm}%

The \ac{CG} method, as presented in Algorithm \ref{alg:cg}, 
is widely regarded as the gold standard iterative solver for sparse \ac{SPD} linear systems
\begin{equation} \label{eq:tridiag}
	\Am \uv = \bv.
\end{equation}
The \ac{CG} method
will produce a sequence of approximations $\uv_k$ for $k=1,2,3,\ldots$
whose residuals we denote by $\rv_k= \bv_k - \Am \uv_k$.
We will rely on the rich a-priori knowledge
on the \ac{CG} method, see e.g.~\cite{HestenesStiefel1952, Axelsson1994, greenbaum1997iterative, meurant1997computation, axelsson2003iteration, saad2003iterative, malek2014preconditioning, greenbaum2021convergence, meurant2024error}.%

In particular, we will use the condition number 
\begin{equation*} 
	\kappa = \kappa({\Am}) = \| \Am \|  \| \Am^{-1} \| = \frac{\lambda_{\max}(\Am)}{\lambda_{\min}(\Am)},
\end{equation*}%
where ${\lambda_{\max}(\Am)}$ 
and     ${\lambda_{\min}(\Am)}$ 
are the largest and the smallest eigenvalue of $\Am$, respectively.
Based on the condition number there is an upper bound on the convergence rate
\begin{equation} \label{eq:conv_bound}
	\frac{|| \rv_k ||}{|| \rv_0 || }  \leq C \; \left( \frac{\sqrt{\kappa}  - 1}{\sqrt{\kappa} +1} \right)^k,
\end{equation}
with $C=\sqrt{\kappa}$.
Additionally, the number of iterations $k$ needed to reduce the error in the energy norm by a factor $\epsilon$ is limited by
\begin{equation} \label{eq:number_of_iterations}
	k \leq \frac12 \sqrt{\kappa(\Am)} \ln (2/ \epsilon).
\end{equation}
Furthermore it is well-known that the \ac{CG} iteration (in exact arithmetic)
will terminate at latest with $k=n$ iterations and then with $\rv_n= 0$.  
We will also make use of the Krylov spaces
\begin{equation} 
	\mathcal{K}_k(\Am, \bv) = \mbox{span} 
		\{
			\bv, \Am \bv, \Am^2 \bv , \ldots , \Am^{k-1}  \bv 
		\}.
\end{equation}%
The \ac{CG} algorithm computes a sequence of approximate solutions $\uv_k \in \mathcal{K}_k(\Am, \bv)$.
As other Krylov subspace methods, it will iteratively 
construct $\uv_k \in \mathcal{K}_k(\Am, \bv)$ such that it satisfies a specific optimality criterion. 
In particular, when the exact solution of the linear system satisfies
$\uv^*= \Am^{-1} \bv \in \mathcal{K}_k(\Am, \bv)$,
the \ac{CG} method will produce this solution exactly (but subject to roundoff errors) in at most $k$ iterations.
This happens at latest when $k=n$, but in practice, 
sufficiently accurate solutions are often found with significantly fewer
iterations $k \ll n$.
Convergence is typically checked by evaluating $|| \rv_k ||$ and stopping 
when it drops below a suitably chosen threshold, as it is 
implemented in line \ref{exit-line} in Algorithm \ref{alg:cg}.

The remainder of this report is organized as follows. 
Section~\ref{sec:intro} introduces a simple illustrative example that exhibits a pronounced \loli effect. 
Section~\ref{sec:other-solvers} extends this discussion through a series of systematic numerical experiments investigating a range of problem settings and iterative solver algorithms. 
Section~\ref{sec:advanced} reports numerical experiments for two-dimensional problems.
Finally, Section~\ref{sec:precond} outlines how an improved performance can be achieved
through preconditioning strategies based on hierarchical transformations that reduce the \loli.
The report concludes with general remarks in Section~\ref{sec:conc}.

\section{A first example illustrating the locality limitation}
\label{sec:intro}
Sparse \ac{SPD} systems frequently arise from the discretization of differential equations.
In these applications, however, the \ac{CG} method may suffer from slow convergence.
The following simple example problem serves to illustrate a situation in which the \ac{CG} method fails to be an efficient solver.
In later sections, we will demonstrate that similar limitations can also affect other iterative methods.

As a starting point, consider the two-point boundary value problem
\begin{equation} \label{eq:example1}
	- u^{\prime\prime} + \gamma^2 u  =  f(x)  \text{~for~}  x \in (0,1)
\end{equation} 
with $\gamma \geq 0$  and subject to the boundary conditions
\begin{equation} \label{eq:bcs}
	u^\prime(0) = 0 \text{~and~ } u(1) = c_r.  
\end{equation}
This is a linear elliptic boundary value problem in 1D with a Neumann boundary condition on the left
and a Dirichlet boundary condition on the right end of the interval, respectively.
We will first choose $f(x)=0$ and $c_r= \cosh{\gamma}$ so that the
unique analytic solution is given by 
\begin{equation} \label{eq:true}
	u^*(x) = \cosh(\gamma x).
\end{equation}
For $\gamma=0$ the solution is particularly simple and becomes $u(x) = 1$.

A standard finite difference discretization on a uniform grid with the $n$ grid points
	$x_i = i h$, for $i= 0,1, 2, \ldots, n$, and $h=1/n$
leads to an \ac{SPD} tridiagonal linear system of the form (\ref{eq:tridiag}) 
where
\begin{equation}
\label{eq:1D-system}
\Am = \Am(h) = 
\frac{1}{h^2} \begin{bmatrix}
d/2 & -1 & 0 & \cdots & 0 & 0 \\
-1 & d & -1 & \ddots & & \vdots \\
0 & -1 & d & \ddots & \ddots & 0 \\
\vdots & \ddots & \ddots & \ddots & -1 & 0 \\
0 & & \ddots & -1& d & -1 \\
0 & 0 & \cdots & 0 & -1 & d
\end{bmatrix} \in \bR^{n \times n} 
\
\text{~and~}
\bv = \bv(h) = \frac{1}{h^2} \begin{bmatrix}
0 \\
0 \\
\vdots \\
\vdots \\
0 \\
0 \\
c_r 
\end{bmatrix} \in \bR^n . 
\end{equation}
The diagonal elements are given by $d= 2+\gamma^2 h^2$.
Note also that the $(n+1)$st discrete solution value $\uv[n]$  
corresponding to the Dirichlet condition at $x_n=1$ 
has already been eliminated from the system.
This reduces the size of the system to dimension $n$ (rather than $n+1$), 
and results in a symmetric system of $n$ unknowns. 
Note also that for $\gamma >0$ the standard Gerschgorin theorem ensures positive definiteness.
In the case $\gamma=0$, all rows have a row sum of $1$, 
except the last one. 
However, also in this case positive definiteness is ensured due to the Dirichlet condition at one interval endpoint.%
Thus the matrix ~(\ref{eq:1D-system}) is suitable to employ the \ac{CG} method as it is presented in Algorithm~\ref{alg:cg}.
\begin{remark} \label{rem:model}
The differential equation~(\ref{eq:example1}) together with the boundary conditions~(\ref{eq:bcs}) represents a prototypical
reaction--diffusion problem that could arise 
in a wide range of physical applications, modeling, e.g., energy and mass transport.
One may, for instance, interpret the solution as the concentration of a chemical species in a porous medium, governed by diffusion combined with a reaction or absorption term. 
In this interpretation, the boundary conditions correspond to a prescribed concentration at the right end of the domain and a zero-flux condition at the left endpoint. 
A natural quantity of interest  in such a setting --and thus also the primary objective of a numerical computation--
could then be the concentration at the left boundary.
\end{remark}
\begin{remark} \label{rem:model2}
Alternatively, the resulting linear system may be viewed as originating from the linearization of a more complex nonlinear model. 
Similar equations also arise in time-dependent problems, 
where the zeroth-order term typically results from the implicit time discretization.
\end{remark}
\begin{remark}
Moreover, systems of this type, then typically in 2D or 3D, are fundamental in computational fluid dynamics,
where they, e.g., appear in pressure-correction schemes that are needed employed to enforce mass conservation.
\end{remark}

For $n=64$ (i.e.~$h=1/64$ ) and $\gamma=2$,
the condition number can be computed numerically as $\kappa(\Am) \approx 2572.46$.
We also recall that here the condition number scales like  
$\kappa(\Am(h)) = {\mathcal O}(h^{-2})$, see, e.g., \cite{leveque2007finite}.
Note also that the bound (\ref{eq:number_of_iterations}) together with the asymptotic behavior of the condition number
predicts that at most $\mathcal{O}(n)$ iterations will be needed to solve system (\ref{eq:tridiag}) to a prescribed accuracy.
This said, we fix the dimension for our example 
to $n=64$ and $h=1/64$
in the remainder of the section.

\begin{remark}
For tridiagonal systems such as~(\ref{eq:tridiag}), 
the most efficient solvers are typically direct methods based on Gaussian elimination, or,
in the \ac{SPD} case, variants of the Cholesky factorization. 
Advanced sparse direct solvers such as in \cite{davis2016survey, duff2017direct} are not necessary for narrowly banded matrices. 
Appropriate techniques include sequential approaches such as the Thomas algorithm~\cite{press2007numerical, conte2017elementary} 
as well as hierarchical and parallelizable schemes based on cyclic reduction~\cite{boisvert1991algorithms, bondeli1994cyclic}. 
Both approaches achieve $\mathcal{O}(n)$ complexity, comparable to that of a single iteration of an iterative method. 
Consequently, iterative solvers are generally not competitive for such problems in terms of either efficiency or attainable accuracy.
\end{remark}
In this work, the simple tridiagonal system~(\ref{eq:tridiag}) is employed as a model problem 
to illustrate and analyze characteristic convergence properties of iterative methods. 
In Sect.~\ref{sec:advanced}, we extend the discussion to two-dimensional example problems and 
demonstrate that the essential observations carry over to this more general setting.
\begin{figure}[th]
    \centering
    \includegraphics[width=0.7\textwidth]{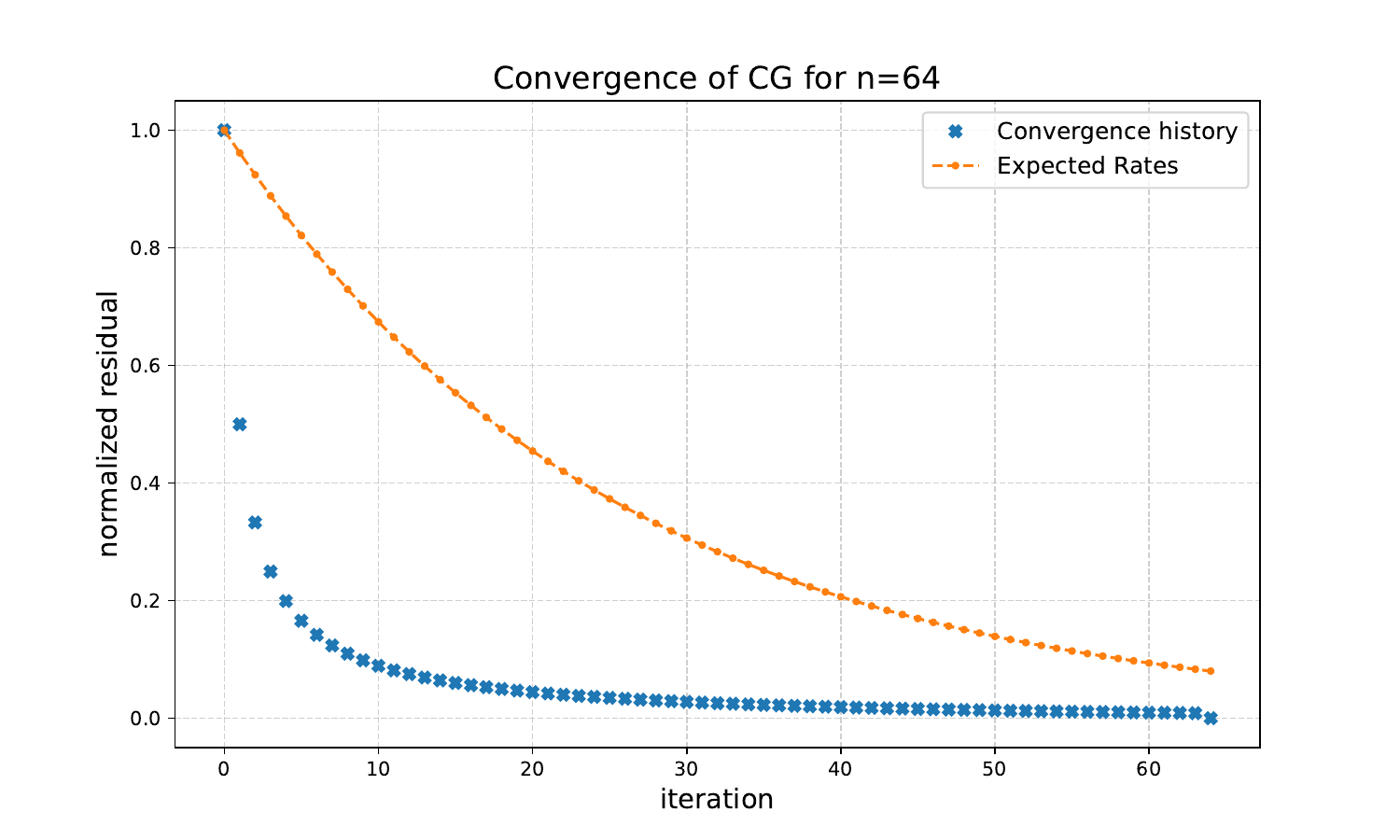}
    \caption{The normalized residual when solving (\ref{eq:example1})
    is displayed depending on the progress of the \ac{CG} iterations. 
    The orange dotted line is the upper bound of the convergence from eq.~(\ref{eq:conv_bound}) with $C=1$.
     }
    \label{fig:convergence_cg_plain}
\end{figure}

For our example with $\gamma=2$, Fig.~\ref{fig:convergence_cg_plain}
shows the convergence of the residual norm towards $0$.
The initial iterations exhibit rapid convergence, followed by a gradual slowdown,
while the residual continues to decrease monotonically, in line with theoretical expectations. 
Using a standard stopping criterion \cite{arioli1992stopping} in line~\ref{exit-line} of Algorithm~\ref{alg:cg}, 
convergence is declared after $k = n = 64$ iterations. 
It is also worth noting that the theoretical convergence bound given in~(\ref{eq:conv_bound}),
included for reference, turns out to be rather pessimistic, as the observed convergence is substantially faster.

At first glance, these results appear entirely satisfactory.
But there are two glitches that should catch our attention.
A careful inspection of Fig.~\ref{fig:convergence_cg_plain} shows a small 
but suspicious kink in the behavior of the residual norm at $k=64$.
Second, despite the initially rapid convergence, the method does not terminate in fewer than $n$ iterations under the default stopping criterion of Algorithm~\ref{alg:cg}.
Becoming curious, we proceed to take a closer look. 

\begin{figure}[th]
    \centering
    \includegraphics[width=0.7\textwidth]{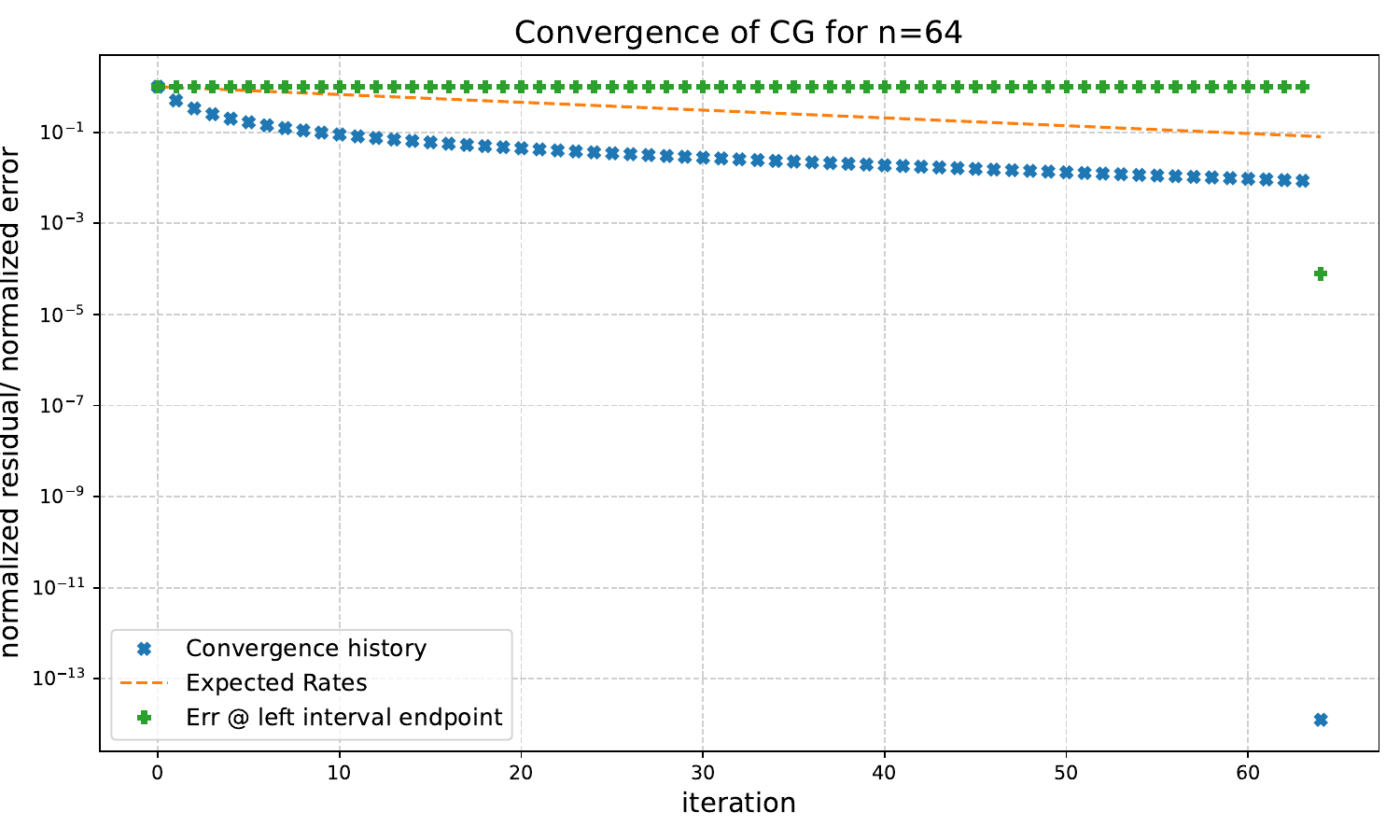}
    \caption{Same as Fig.~\ref{fig:convergence_cg_plain} showing the reduction of the normalized residual 
   for example (\ref{eq:example1}) with $f(x)=0$, now in semilogarithmic scale.
    Note the sudden drop to $|| \rv_n || \leq 10^{-13}$ occurring in the last iteration.
    Additionally, the error at the left endpoint is displayed (marked by green + symbols).}
    \label{fig:convergence_cg_logarithmic}
\end{figure}

To this end, we turn to a semilogarithmic representation. 
The same data as in Fig.~\ref{fig:convergence_cg_plain}.
are also displayed in Fig.~\ref{fig:convergence_cg_logarithmic}.
Presented in this (better!) way, it becomes obvious that there is a problem.
When viewed on a logarithmic scale, the convergence behavior appears markedly different.
For all iterations $k < n$, the residual norm exhibits near-stagnation, remaining above $10^{-3}$, 
before dropping abruptly to approximately $10^{-13}$ at $k = n$. 
Thus, the substantial reduction of the residual occurs almost exclusively in the final iteration.
Fig.~\ref{fig:convergence_cg_logarithmic} also reports the evolution of the error 
in the iterates relative to the analytical solution of~(\ref{eq:example1}), 
evaluated at the first component, corresponding to $x=0$. 
Interestingly, this pointwise error, $\ev_k[0] := u^*(0) - \uv_k[0]$, 
remains constant at a value of $1$ throughout all intermediate iterations and decreases only in the final step, where it reaches approximately $10^{-4}$. 
Since this quantity is measured with respect to the analytical solution, 
the remaining error after convergence of the \ac{CG} method reflects the discretization error at $x = 0$.

To better understand what is happening, we proceed further and additionally visualize the intermediate iterates.%
\begin{figure}[th]
    \centering
    \includegraphics[width=0.7\textwidth]{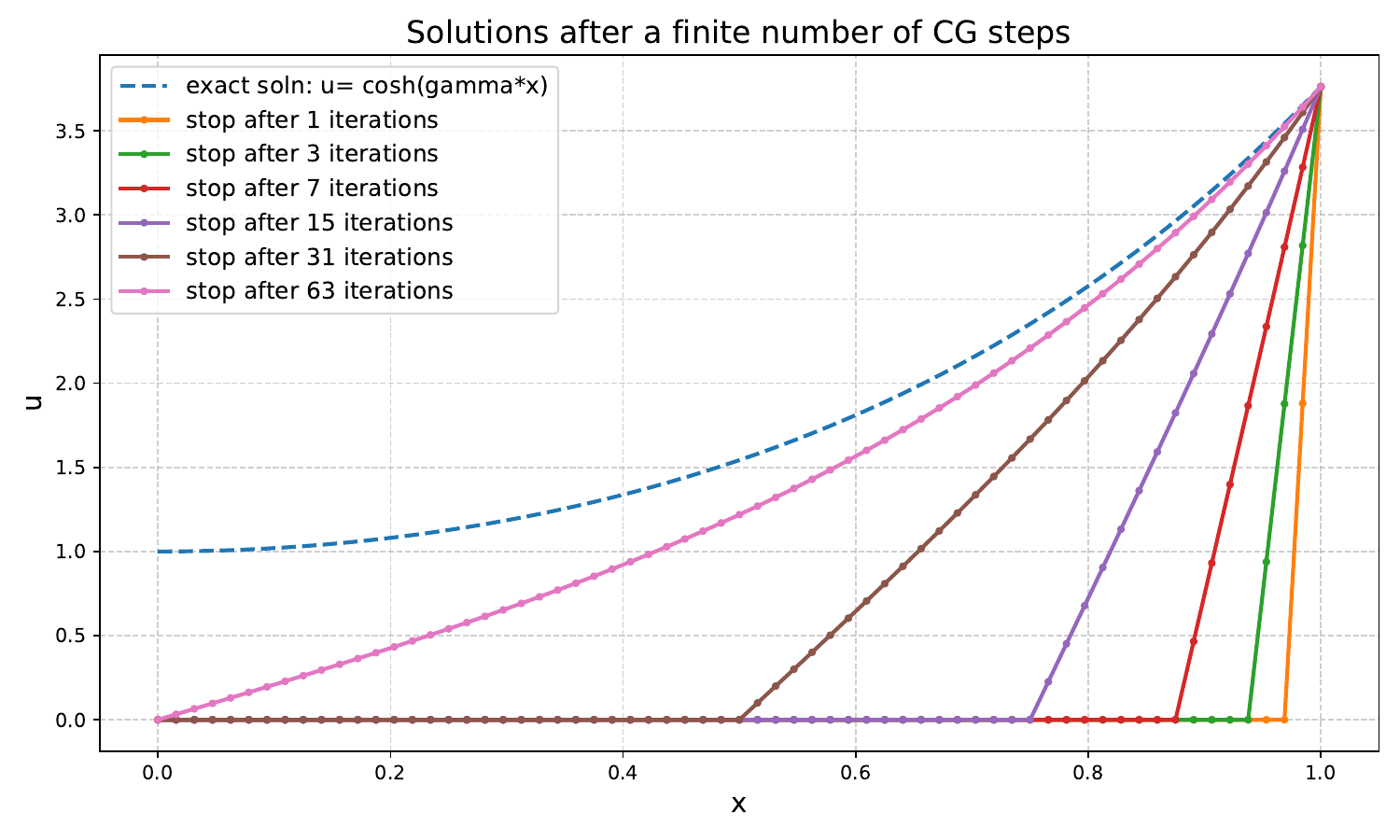}
    \caption{The true solution compared with the results after $k= 1, 3, 7, 15, 31, 63$ \ac{CG} iterations.}
    \label{fig:stopping_64}
\end{figure}

In Fig.~\ref{fig:stopping_64}, we observe that all iterates with $k < n$ not only exhibit relatively large residuals, 
but also remain qualitatively incorrect. 
In particular, for all intermediate iterations the solution satisfies $\uv_k[0] = 0$ at the left endpoint of the interval. 
Consequently, the approximation at $x = 0$ does not improve toward the expected value 
$u^*(0) = \cosh(0) = 1$ until the final iteration is reached.
A more detailed inspection shows that each iterate $\uv_k$ satisfies $\uv_k[j] = 0$ for all $j < k$, 
corresponding to the inaccurate approximation $u(x) = 0$ on the interval $x \in (0, 1 - kh)$. 
Thus, in the $k$-th iteration the approximate solution remains identically zero over a subinterval 
whose length decreases only gradually, shrinking by a single grid spacing $h$ with each iteration.
In contrast, the exact solution satisfies the fundamental property $u(x) > 0$ throughout the domain. 
This information enters the system at the right boundary, 
where the Dirichlet condition enforces $u(1) > 0$,
but this feature propagates only slowly into the interior.
Indeed, we observe that in each iteration of the \ac{CG} method
 this information advances by just one grid point from right to left. 
 Only in the final iteration, when $k = n$, does it reach the left boundary at $x = 0$,
 and only then does the method produce a qualitatively correct approximation 
 to the solution $u(x) = \cosh(\gamma x)$.%

Several conclusions can already be drawn from this simple illustrative example.
\begin{itemize}
\item
The \ac{CG} method exhibits a pronounced limitation in the speed at which information propagates through the computational domain.
\item
	The sparsity pattern of the system matrix $\Am$ imposes a fundamental upper
	bound on the rate of information transfer within the algorithm.
\item
	In each iteration of the \ac{CG} method, information propagates by at most one mesh point;
	consequently, $n$ iterations are required before a Dirichlet condition imposed at the right boundary
	can influence the solution at the left endpoint.
\item
	A qualitatively accurate approximation of the solution cannot be obtained before the final iteration,
	that is, when $k = \dim(\mathbb{R}^n)$.
\item
	This behavior reflects an inherent \loli of the \ac{CG} algorithm.
\end{itemize}
The resulting slow convergence depends directly on the number of grid points $n$ and, equivalently, 
on the mesh size $h$ used to discretize~(\ref{eq:example1}). 
As the mesh is refined to improve discretization accuracy, 
the number of iterations required for convergence increases accordingly, 
further delaying the attainment of acceptable solution accuracy.
\begin{remark}
The \emph{locality limitation} of the \ac{CG} method is well known and has been studied in particular for graph Laplacians; 
see, for example,~\cite{SpielmanCGDiameter, hu2024solving}. 
In this setting, the close relationship between the graph diameter, the spectrum of the Laplacian, 
and the resulting condition number is used to motivate the necessity of preconditioning when solving
problems associated with graphs of large diameter. 
Somewhat surprisingly, however, there appears to be no systematic investigation of this phenomenon
for more general linear systems, such as those arising from the discretization of differential equations.
\end{remark}

Thus the introductory example naturally raises several questions:
\begin{itemize}
\item
	Is the observed behavior specific to this constructed example, or does it also arise in practically relevant applications?
\item
	How does this effect manifest itself in two- and three-dimensional problems,
	where the \ac{CG} method is commonly employed?
\item
	To what extent can the slow convergence be alleviated through the use of suitable preconditioners?
\end{itemize}
These questions will be addressed in the remainder of the paper.

\section{A broader experimental study of the \loli}
\label{sec:other-solvers}
In this section we will study whether the observed behavior is representative for iterative methods more generally
or whether we have only found a singular example in which the iteration is untypically slow.

\subsection*{More general right hand sides}
Example~(\ref{eq:example1}) with $f(x)=0$ gives rise to a linear system 
whose right-hand side vector $\bv$ is proportional to the last unit vector in~(\ref{eq:1D-system}). 
Consequently, the exact solution
	$\uv^* = \Am^{-1}\bv$
is a scaled version of the final column of $\Am^{-1}$ and 
can therefore be interpreted as a discrete representation of the associated Green’s function.
To explore whether this alone leads to the slow convergence, we modify the problem by instead choosing $f(x)=10$. 
This results in a right-hand side vector $\bv$ with strictly positive entries in all components, so that the right hand side is
not strictly localized any more.
The corresponding analytical solution of~(\ref{eq:example1}) is then given by $u(x)=10+\cosh(\gamma x)$, 
when the Dirichlet boundary condition is adjusted accordingly to $u(1) = c_r = 10+\cosh(\gamma)$.
\begin{figure}[ht]
    \centering
    \includegraphics[width=0.7\textwidth]{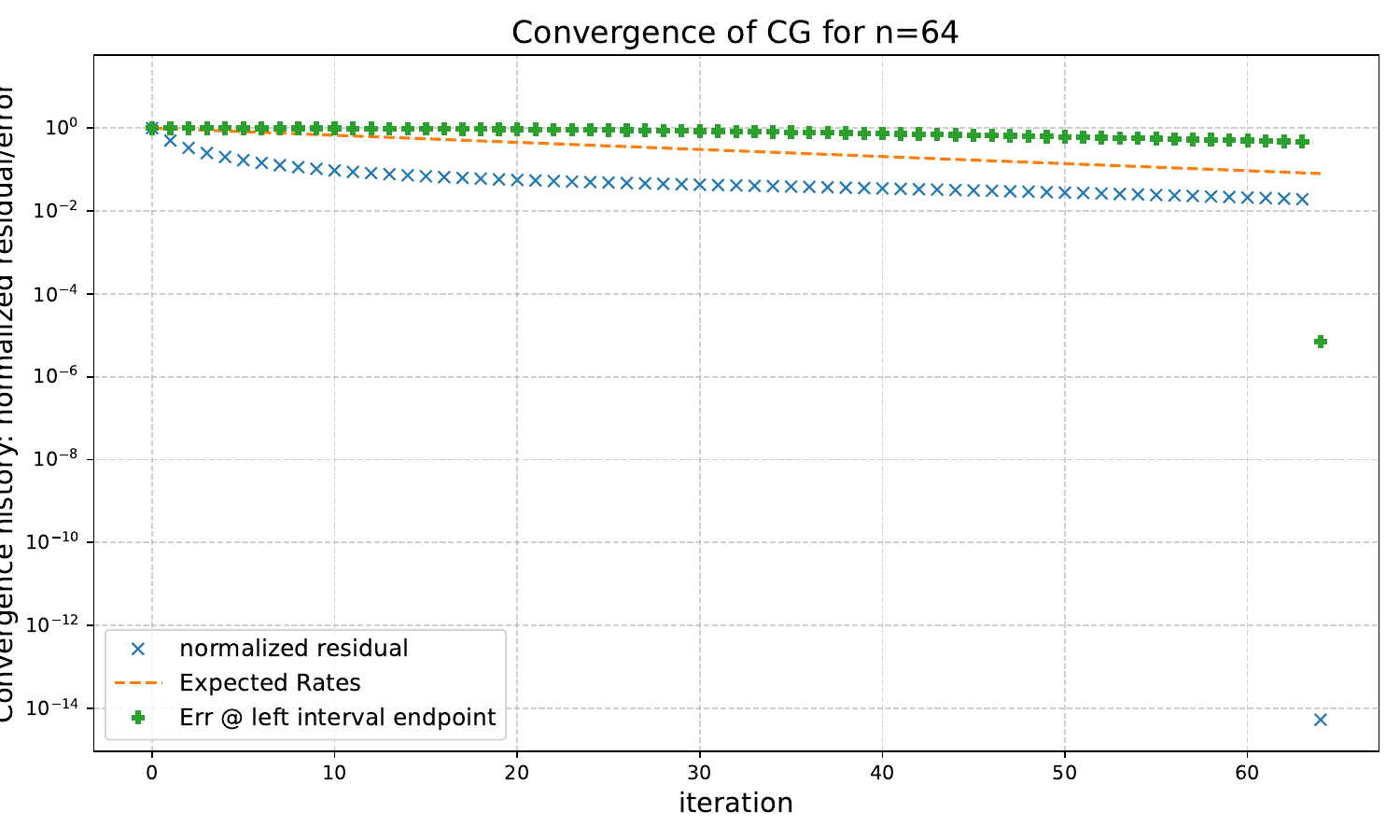}
    \caption{The reduction of normalized residual for example (\ref{eq:example1}) with $f(x)=10$ solved with \ac{CG} iterations.
    Additionally the error in the left endpoint is shown. 
    Note again the sudden drop of the residual norm and also the error when $k=n$.
    \label{fig:exa2-convergence}}
\end{figure}
\begin{figure}[ht]
    \centering
    \includegraphics[width=0.7\textwidth]{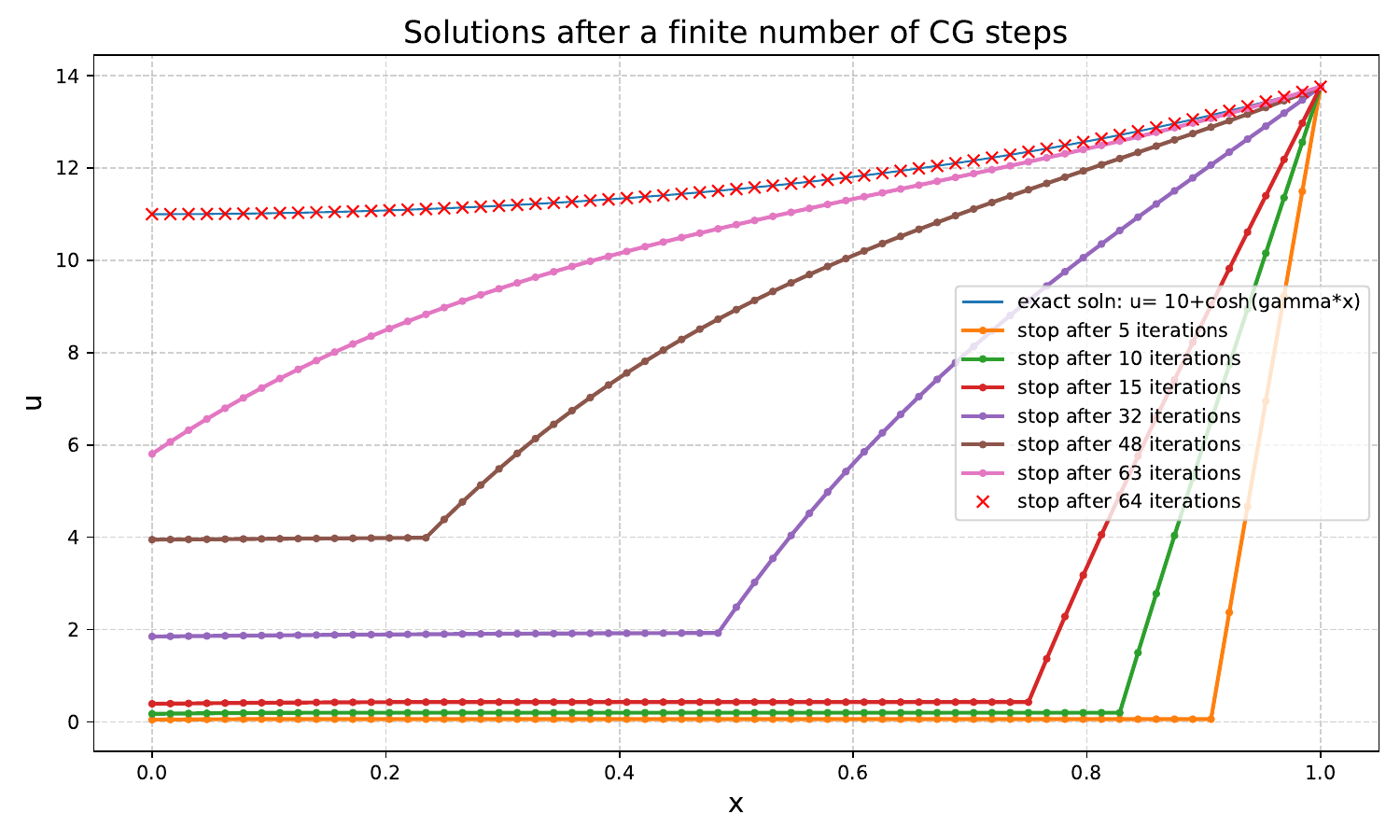}
    \caption{The true solution compared with the intermediate results after $k= 5, 10, 15, 32, 48, 63, 64$ \ac{CG} iterations.
    \label{fig:exa2-stopping}}
\end{figure}

The corresponding convergence behavior and selected intermediate iterates are shown in
Figs.~\ref{fig:exa2-convergence} and~\ref{fig:exa2-stopping}, respectively.
The residual convergence exhibits a pattern similar to that observed previously, characterized by an extended phase of stagnation.
A substantial reduction of the residual again occurs only in the final iteration, when $k = n$ is reached. 

When viewed on a logarithmic scale, as in the figures, the onset of convergence warrants a careful inspection.
In contrast to the previous example, the associated Krylov spaces are less restrictive. 
They allow that the left endpoint assume values $\neq 0$ already before the final iteration. 
Accordingly, unlike in Fig.~\ref{fig:convergence_cg_logarithmic}, 
the intermediate iterates now satisfy $\uv_k[0] > 0$ at the left boundary. 
Nevertheless, although these values increase monotonically toward the correct solution, 
their convergence remains very slow.

To further illustrate the behavior of the \ac{CG} algorithm, 
selected intermediate approximations are shown in Fig.~\ref{fig:exa2-stopping}.
This sequence clearly illustrates the monotonic convergence of the solution at the left endpoint of the interval. 
While this behavior marks some improvement compared to Fig.~\ref{fig:stopping_64}, 
the overall effect is less pronounced than one might expect. 
The convergence remains sufficiently slow that all iterates with $k < n$ would still be regarded as qualitatively inadequate.

\subsubsection*{Better initial guesses}
\begin{figure}[ht]
    \centering
    \includegraphics[width=0.7\textwidth]{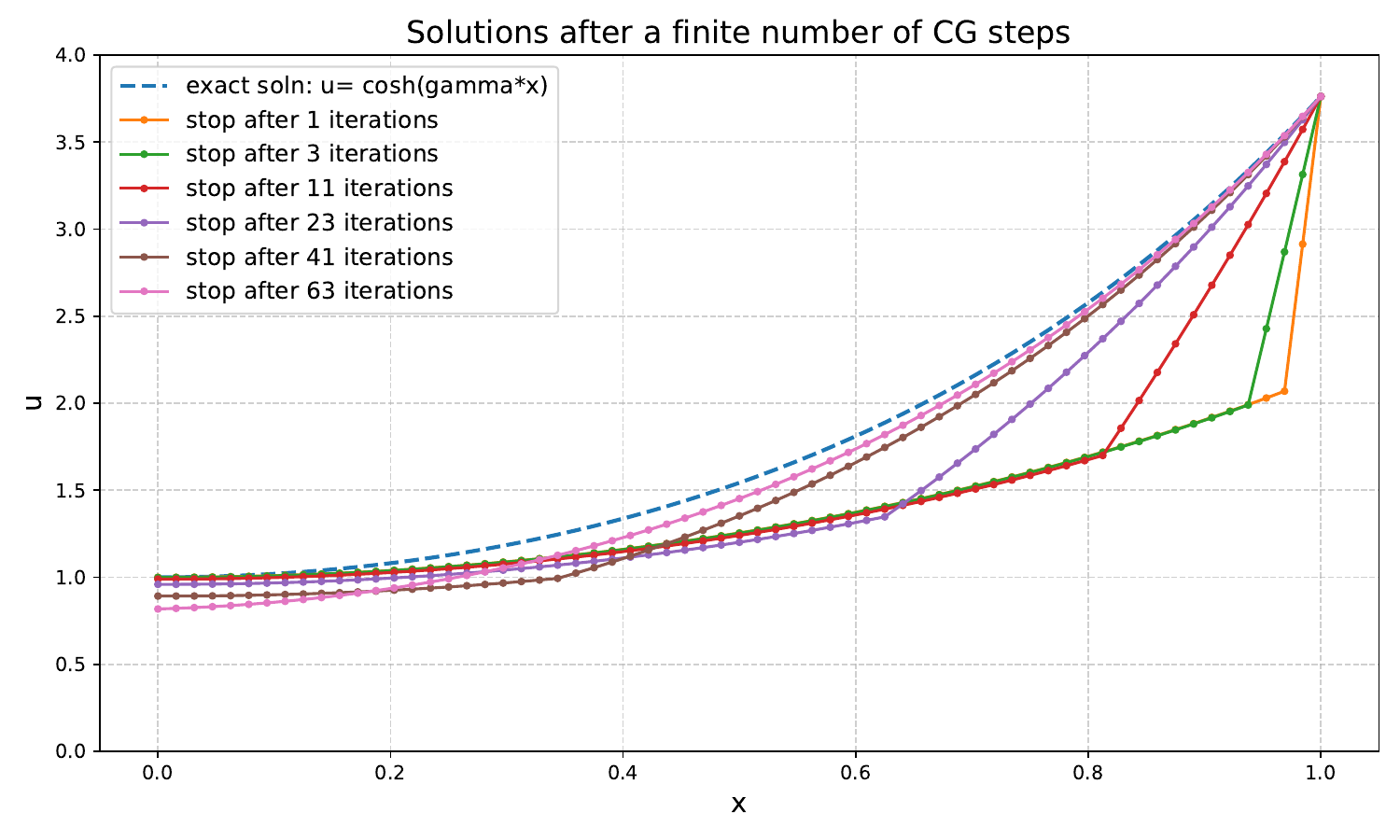}
    \caption{Intermediate iterates for the \ac{CG} iteration applied to (\ref{eq:tridiag}) 
    but starting with initial guess $\uv_0$ representing $\cosh(0.7 \cdot x)$.
    Note that at the left end of the interval ($x=0$), the initial guess has the correct value $\uv_0[0] = 1.0$, 
    but that the iterates are deviating from this value.
    In fact we observe that the last iterate before convergence of the \ac{CG} iteration, i.e.~$\uv_{63}$, exhibits the worst error at $x=0$.} 
    \label{fig:stopping-with-initial-guess}
\end{figure}

In the previous examples, the \ac{CG} iteration was initialized with the canonical choice $\uv_0 = 0$.
In many practical situations, however, more informed initial guesses are available. 
This is commonly the case when the solution of a linear system is embedded in an outer iterative process,
such as a nonlinear iteration or a time-stepping scheme.

To examine the influence of the initial guess on convergence, we consider here a modified initialization of the \ac{CG} algorithm.
While a systematic exploration of all possible initial vectors is beyond the scope of this work, 
the following example provides some illustrative insight. 
Specifically, we choose $\uv_0$ as a discrete representation of~(\ref{eq:true}) with $\gamma = 0.7$,
reflecting a scenario in which a sequence of problems of the form~(\ref{eq:example1}) is solved for increasing values of $\gamma$.
The resulting iterates are shown in Fig.~\ref{fig:stopping-with-initial-guess}.

Notably, here the initial guess attains the correct value at the left boundary, satisfying $\uv_0[0] = u(0) = 1$.
However, this property is not preserved during the iteration. 
Instead, and somewhat counterintuitively, the intermediate iterates exhibit increasingly
large errors at $x = 0$ until convergence is eventually achieved in the $n$-th iteration.
Only then, the \loli has been overcome. 
For this example, we omit the corresponding convergence histories for brevity.

\subsection*{Matrices with better conditioning}
The behavior discussed so far corresponds to rather poorly conditioned systems, with
$\kappa(\Am) = \mathcal{O}(n^2)$. 
For the examples considered above with $n=64$, we obtain
numerically $\kappa(\Am) \approx 2572$, for which the theory already predicts slow convergence,
for instance as indicated by the bound~(\ref{eq:conv_bound}).
\begin{remark} 
By increasing the parameter $\gamma$ in~(\ref{eq:example1}), one may construct examples with
significantly improved conditioning. Indeed, values such as $\gamma^2 = \mathcal{O}(h^{-1})$
or even $\gamma^2 = \mathcal{O}(h^{-2})$ naturally arise, for example, when the linear system
results from an implicit time step in a time-dependent problem. In such regimes, faster
convergence is expected.
\end{remark}%

Here, we return to the canonical initial guess $\uv_0 = 0$ and choose $\gamma = 8$,
which yields a condition number of approximately $\kappa(\Am) \approx 252$.
For this case, the theoretical upper bound on the per-iteration residual reduction factor
given by~(\ref{eq:conv_bound}) is approximately $0.88$.
\begin{figure}[ht]
    \centering
    \includegraphics[width=0.7\textwidth]{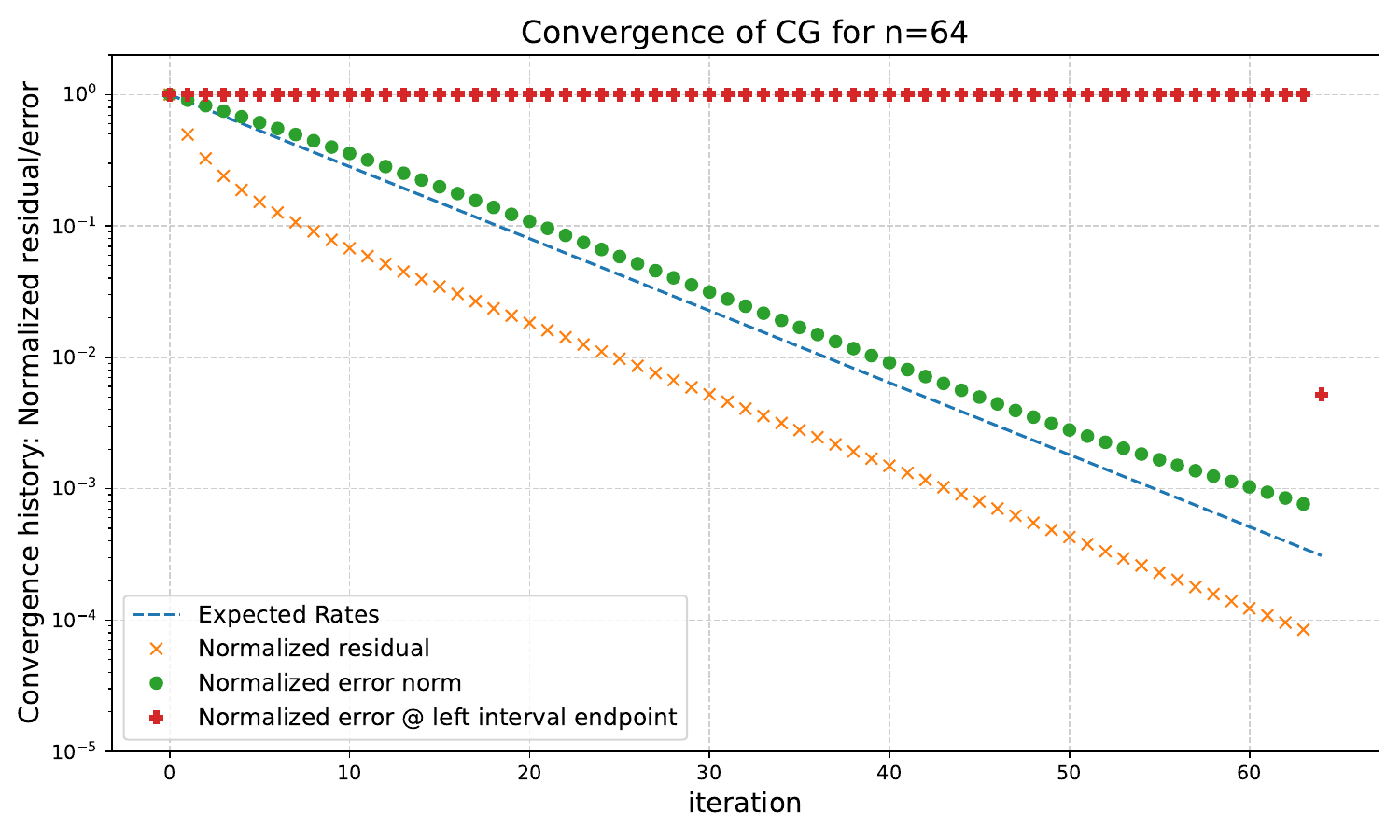}
    \caption{Convergence history of normalized residual norm and normalized error norm 
    in case of $\gamma=8$ and $\kappa(\Am) \approx 252$.
    Note that both values drop to values around $10^{-14}$ at the 64-th iteration, outside the displayed scale.
    The bound on the convergence is displayed for comparison as a dashed line.
    For comparison, also the point-wise error at the left endpoint is shown.} 
    \label{fig:better_conditioned_convergence_history}
\end{figure}
Figure~\ref{fig:better_conditioned_convergence_history} confirms the anticipated acceleration of
convergence. 
The asymptotic rate is accurately predicted by the reduction factor given
in~(\ref{eq:conv_bound}), as indicated in the figure. 
In addition, we here also report the normalized
Euclidean error norm of the iterates relative to the reference numerical solution; this
quantity decays at a comparable rate.

In contrast, the pointwise error at the left endpoint of the interval displays a markedly
different behavior. 
Despite the improved conditioning, it remains essentially unchanged
from the previously observed cases with higher condition number. 
Again, the solution value at the left endpoint stagnates at $0$ throughout the iterations and improves only in the final step,
where it jumps to a value close to $1$. 
The remaining discrepancy at that point is then
dominated by the discretization error.
\begin{figure}[ht]
    \centering
    \includegraphics[width=0.7\textwidth]{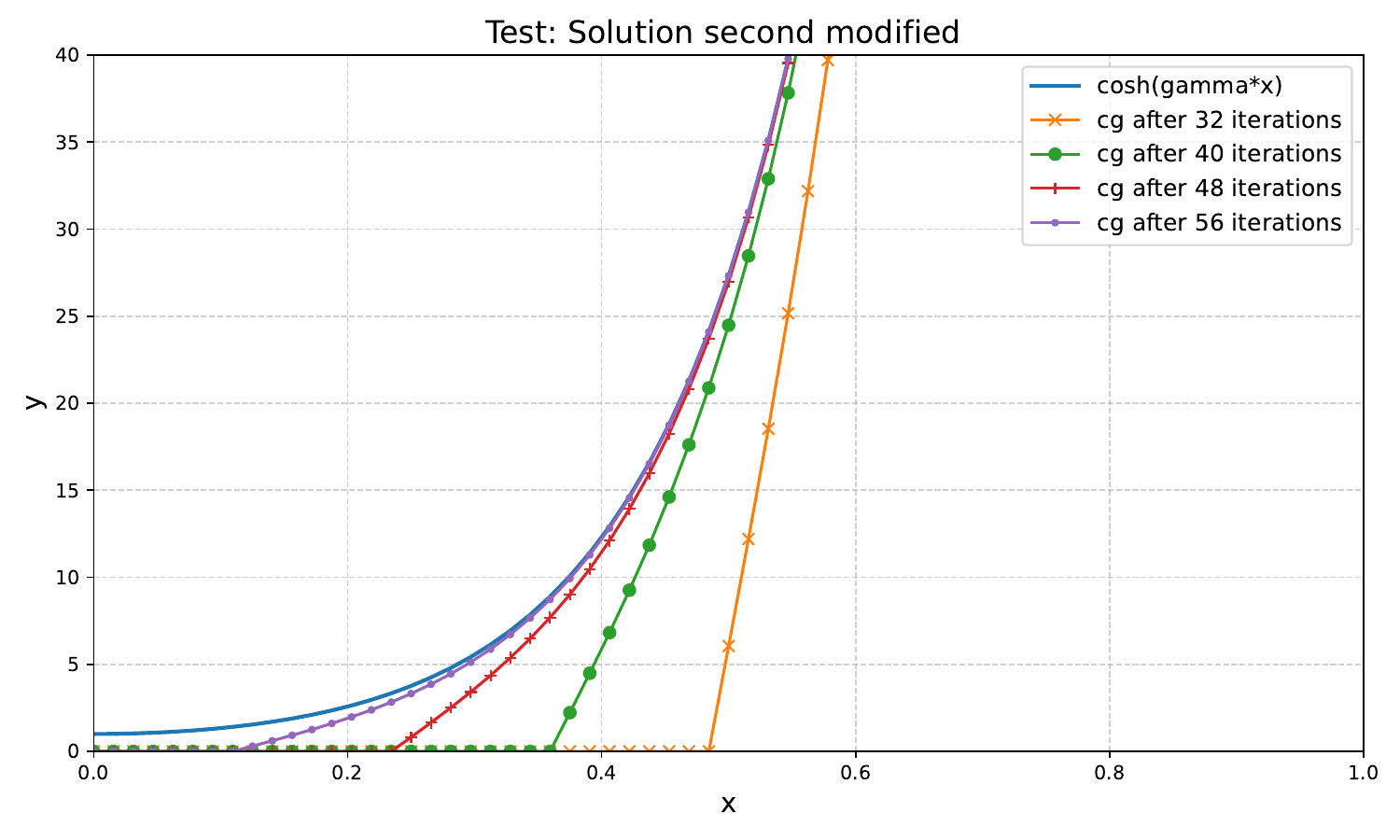}
    \caption{Intermediate solutions obtained by the \ac{CG} algroithm in the case $\gamma=8$ and $\kappa(\Am) \approx 252$.
    For $x \rightarrow 1$ the solutions are outside the displayed scale.} 
    \label{fig:fig:better_conditioned_intermediate}
\end{figure}%

To further illustrate this case, we again visualize selected intermediate iterates produced by
the \ac{CG} algorithm in Fig.~\ref{fig:fig:better_conditioned_intermediate}.
The \emph{locality limitation} remains clearly visible: information imposed by the Dirichlet
boundary condition propagates into the interior only one mesh point per iteration.
A notable difference, however, is that the boundary value $u(1)=\cosh(\gamma)$ is now large and
lies outside the range displayed in Fig.~\ref{fig:fig:better_conditioned_intermediate}.
As a consequence, the distribution of error and residual is no longer balanced across the domain.
Both quantities are dominated by contributions from the right part of the interval, where the
solution magnitude is largest.
During the \ac{CG} iteration, this dominant error component is reduced first, which largely
explains the improved global convergence rates observed in the residual norm.
Note also that Fig.~\ref{fig:better_conditioned_convergence_history} shows the normalized error and residual norms.
In contrast, the behavior in the left part of the domain remains essentially unchanged from the
previous examples: there, the error is barely reduced until the final iteration $k=n$ is reached.
In the global norm, this deficit of the iterates is hidden behind  the larger error components.

\begin{remark}
Whether this behavior constitutes a serious limitation depends strongly on the application and on
the specific computational objective.
If the primary goal is the reduction of a global error measure, such as the residual or the
Euclidean error norm, then the improved condition number indeed yields a reasonably rapid
convergence.
In this case, a solution obtained after, for example, $k \approx n/2$ iterations might already be
judged satisfactory.

If, however, the objective is the accurate evaluation of a local quantity, such as the value
$u(0)$ (see Remark~\ref{rem:model}), then the \loli remains as critical as before,
despite the improved conditioning of the matrix.
\end{remark}
\begin{remark}
In practical applications, premature termination can have significant consequences.
For example, the predicted onset or breakthrough of a pollutant—represented in our model by the value $u(0)$—may be severely underestimated in simulations of transport through porous media.
\end{remark}
In the remainder of this section, we will turn to alternative iterative methods and explore
to which extent they also suffer from a \loli.

\subsection*{Stationary iterative solvers}
\begin{figure}[ht]
    \centering
    \includegraphics[width=0.7\textwidth]{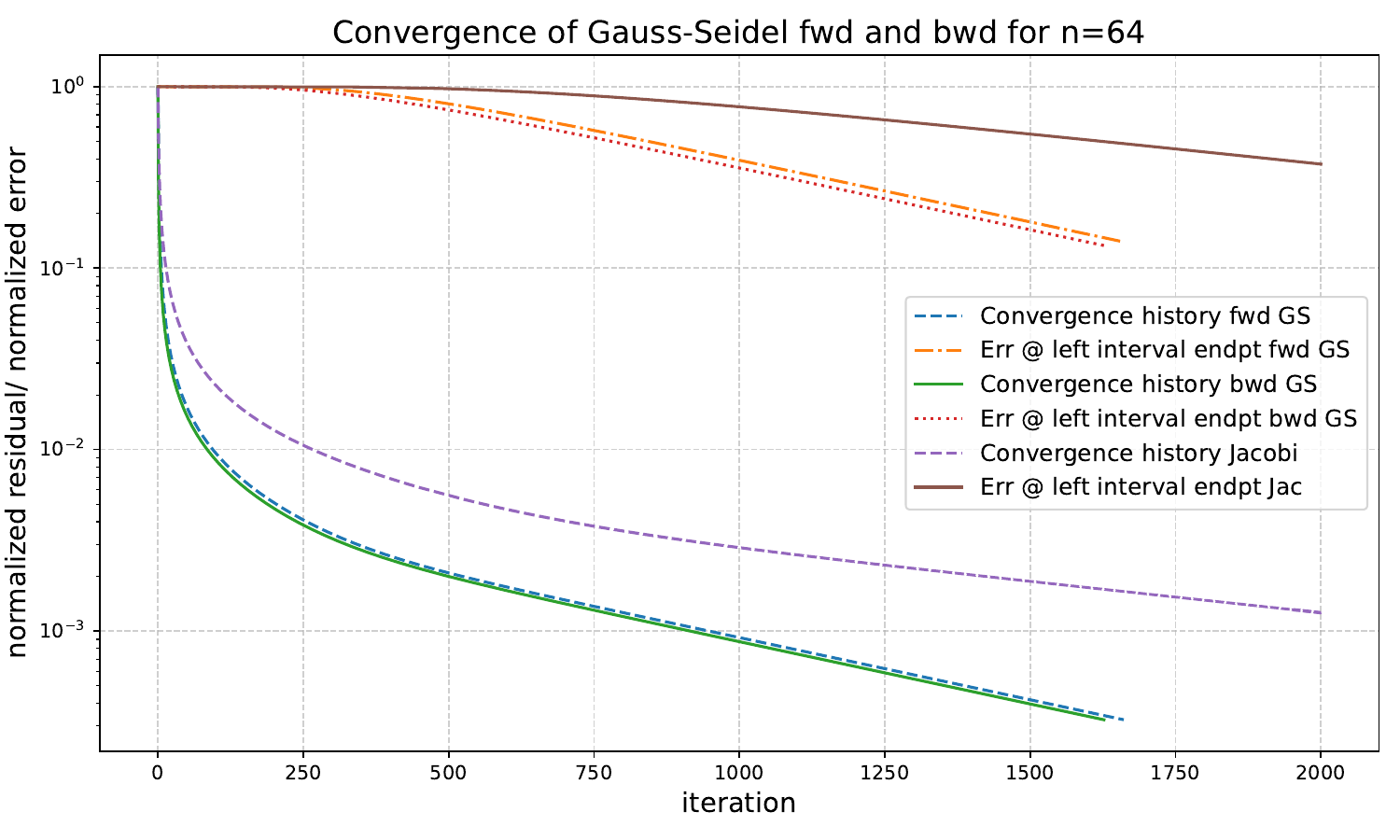}
    \caption{Convergence history for forward and backward Gauss-Seidel- as well as the Jacobi-iteration.%
    \label{fig:gs-rates}} 
\end{figure}

\begin{figure}[ht]
    \centering
    \includegraphics[width=0.7\textwidth]{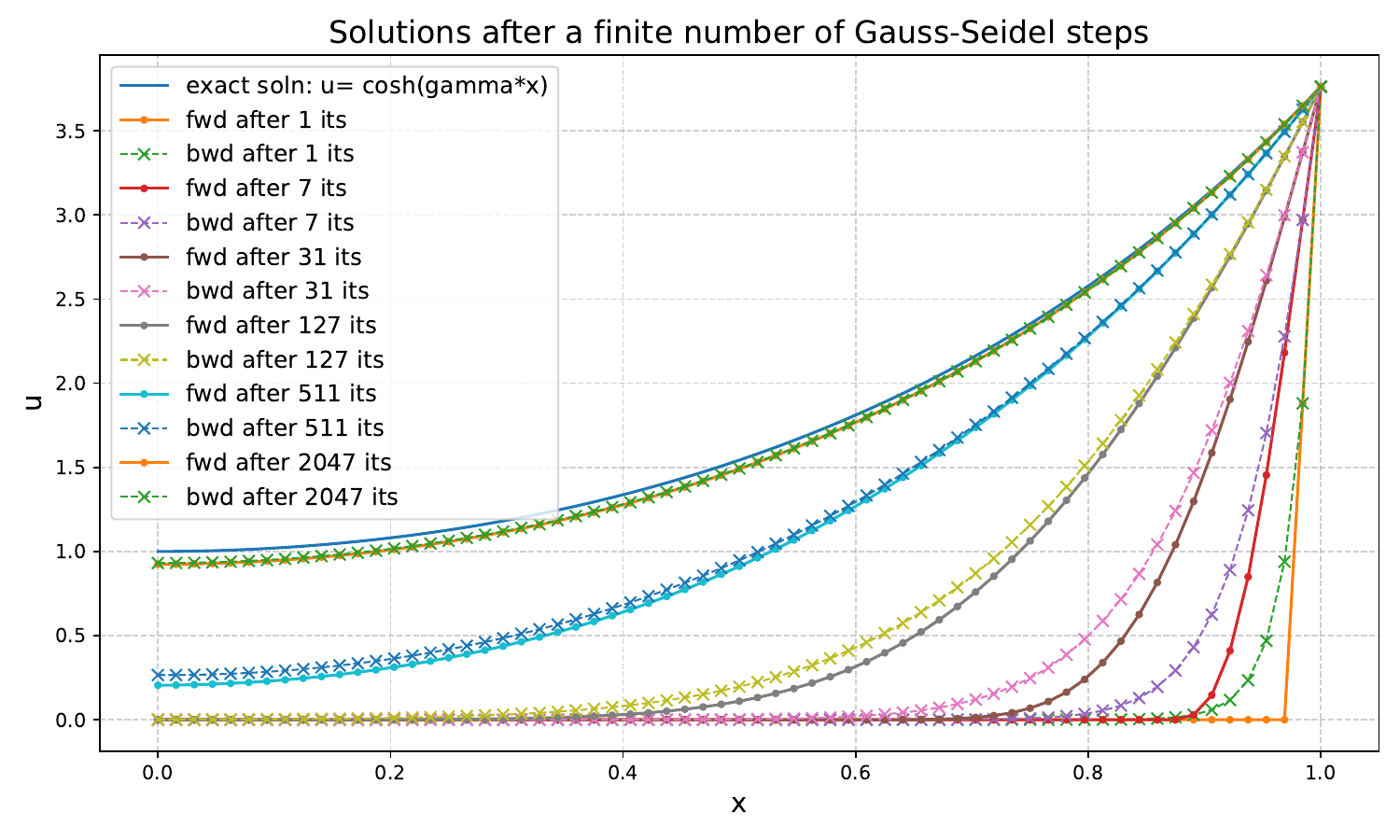}
    \caption{Intermediate solution of the Gauss-Seidel iterations
    when the iteration is stopped as indicated by the legend of the diagram.
    \label{fig:gs-stopping}}
\end{figure}

For a start, we will study stationary iterative methods of the form
\begin{equation} \label{eq:stationary-iteration}
	\uv_{k+1}= \uv_k + \Bm ( \bv - \Am  \uv_k)  \text{~for~} k=0, 1, 2, 3, \ldots ,
\end{equation}
where $\Bm$ is an approximate inverse of $\Am$.
For the relaxed Jacobi method $\Bm = \rho \Dm^{-1}$ where $\Dm$ is the diagonal of $\Am$ and $\rho >0 $ is a relaxation factor.
Similarly, the forward and backward Gauss-Seidel methods can be defined using
$\Bm= (\Dm+ \Lm)^{-1}$ or $\Bm = (\Dm+ \Um)^{-1}$, respectively.
Here we use the splitting of $\Am= \Lm + \Dm + \Um$, where $\Lm$ and $\Um$ 
are the lower and upper triangular sub-matrices of $\Am$, respectively. 
We evaluate the convergence behavior by an experimental setup equal to our evaluation of the \ac{CG} method.
Figs.~\ref{fig:gs-stopping} and \ref{fig:gs-rates} display the results for $n=64$, where we have
returned to the original example of eq.~(\ref{eq:example1}) with $f(x)=0$.
All three methods converge considerably more slowly than the \ac{CG} method, requiring over
$2000$ iterations to reach a qualitatively acceptable solution.

In contrast, the normalized residual first decays quickly over approximately 100 iterations,
but then the reduction also reaches an asymptotic limit.
Forward and backward Gauss-Seidel show very similar performance, while the
Jacobi iteration, as expected, converges at roughly half the rate of Gauss-Seidel.

As before, the error at $x=0$ for all three methods stagnates in the beginning of the iteration.
This is clearly again caused by the \loli. 
Only when it has been overcome, a slow convergence sets in. 
However, even after $2000$ iterations, the residual has been reduced by only three to four orders of magnitude,
which is still insufficient to achieve the discretization-level accuracy
(cf. Fig.~\ref{fig:convergence_cg_logarithmic}).

Next, Fig.~\ref{fig:gs-stopping} presents selected intermediate iterates of the Gauss--Seidel
method for both forward and backward orderings. 
Even for the largest iteration count shown, $k=2047$, 
the computed solutions still differ visibly from the reference solution, which is
included for comparison. 

\begin{figure}[ht]
    \centering
    \includegraphics[width=0.7\textwidth]{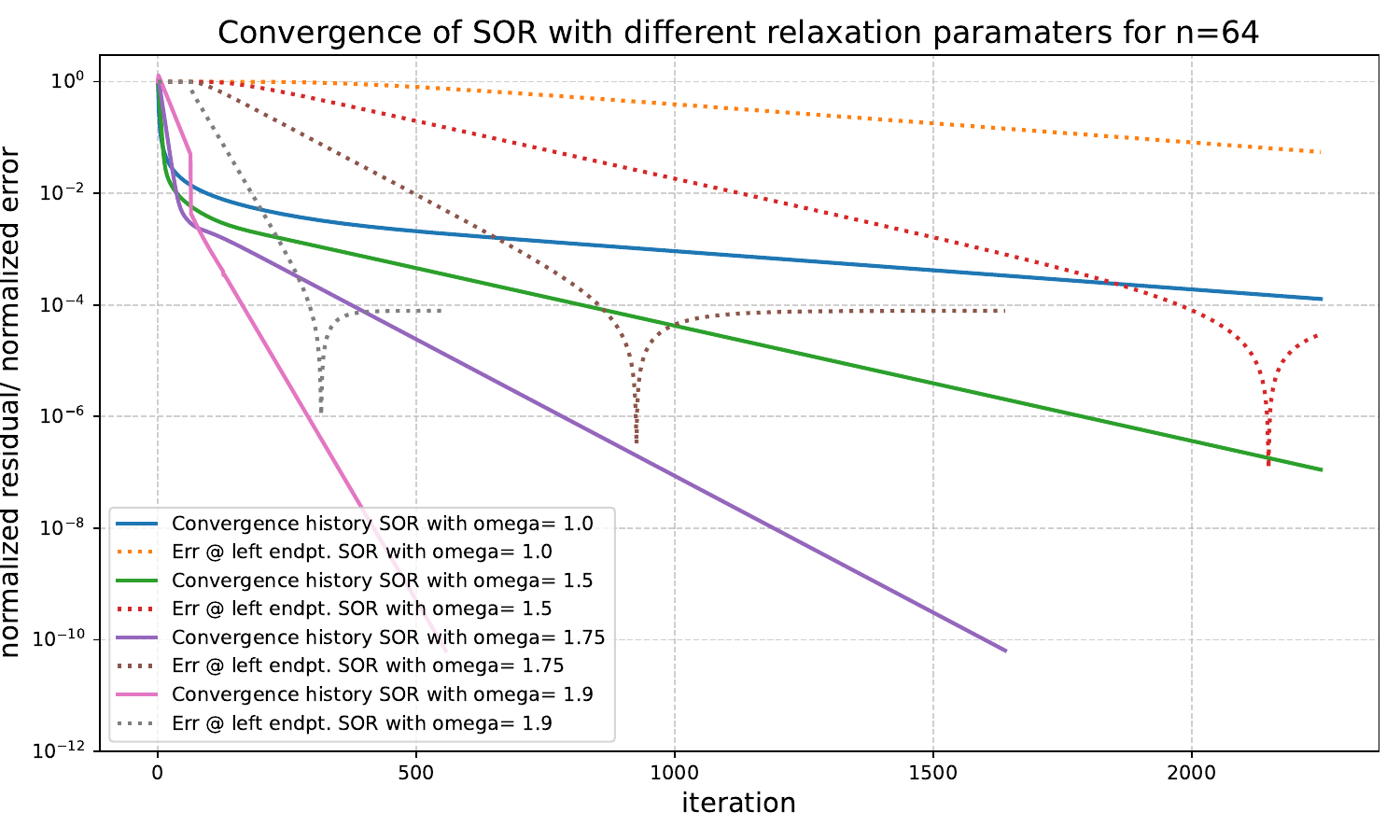}
    \caption{Convergence history of the SOR method for different relaxation parameters
    and 
    for the normalized residual norm as well as the error at the left endpoint.    
    \label{fig:sor}}
\end{figure}
A simple and effective acceleration of the Gauss-Seidel iteration is the classical \ac{SOR} method.
The optimal parameter for \ac{SOR}
may be difficult to find in realistic computations.
Here, we display the convergence history in Fig.~\ref{fig:sor} for different relaxation parameters.
The case $\omega=1.0$ corresponds to the standard Gauss--Seidel method, for which no
acceptable approximation is obtained within $2500$ iterations.

For larger values $\omega>1$, substantially faster convergence is observed.
Nevertheless, independent of the chosen relaxation parameter, the error at the left
endpoint exhibits the same characteristic stagnation up to $k=64$ iterations, which is
a direct consequence of the \loli.
\begin{remark}
An initially counterintuitive behavior is the pronounced kink in the error curve, where
the error temporarily decreases below its final level before converging asymptotically
to approximately $10^{-4}$, corresponding to the discretization error.
This behavior is explained by a sign change of the pointwise error during the iteration:
while the early iterates underestimate the analytical value, 
the final discrete solution satisfies $\uv_n[0] > u(0) = 1$.
The resulting non-monotonicity is therefore an artifact caused by the point wise evaluation of the error.
\end{remark}
\begin{remark}
For the present example with $n=64$ and $\gamma=2$, a nearly optimal \ac{SOR} parameter
is $\omega \approx 1.9$.
With this choice, the \ac{SOR} method reaches discretization-error accuracy in fewer than
$500$ iterations.
We recall that \ac{SOR} with the optimal parameter attains the
same asymptotic complexity as the \ac{CG} method \cite{axelsson2001finite}.
Since stationary iterative schemes are not the primary focus of this work, we do not
pursue this further.
\end{remark}

\subsection*{The \ac{GMRES} method with restarts}
The final method to be discussed in this section is again of Krylov type,
the \ac{GMRES} method \cite{saad1986gmres, greenbaum1996any}.
\ac{GMRES} is arguably the most popular Krylov space method for non-symmetric matrix equations.
We will here employ the method to solve 
the linear system (\ref{eq:1D-system}), as before.
\begin{remark}
In practice,  it is not advised to use \ac{GMRES} for \ac{SPD} systems, since the \ac{CG} method
is considered a better alternative in these cases.
Nevertheless, for our purpose of illustration, the performance and convergence behavior of \ac{GMRES}
for the \ac{SPD} example
is interesting, since we can compare the results for the same example problem.
\end{remark}
\begin{figure}[th]
    \centering
    \includegraphics[width=0.7\textwidth]{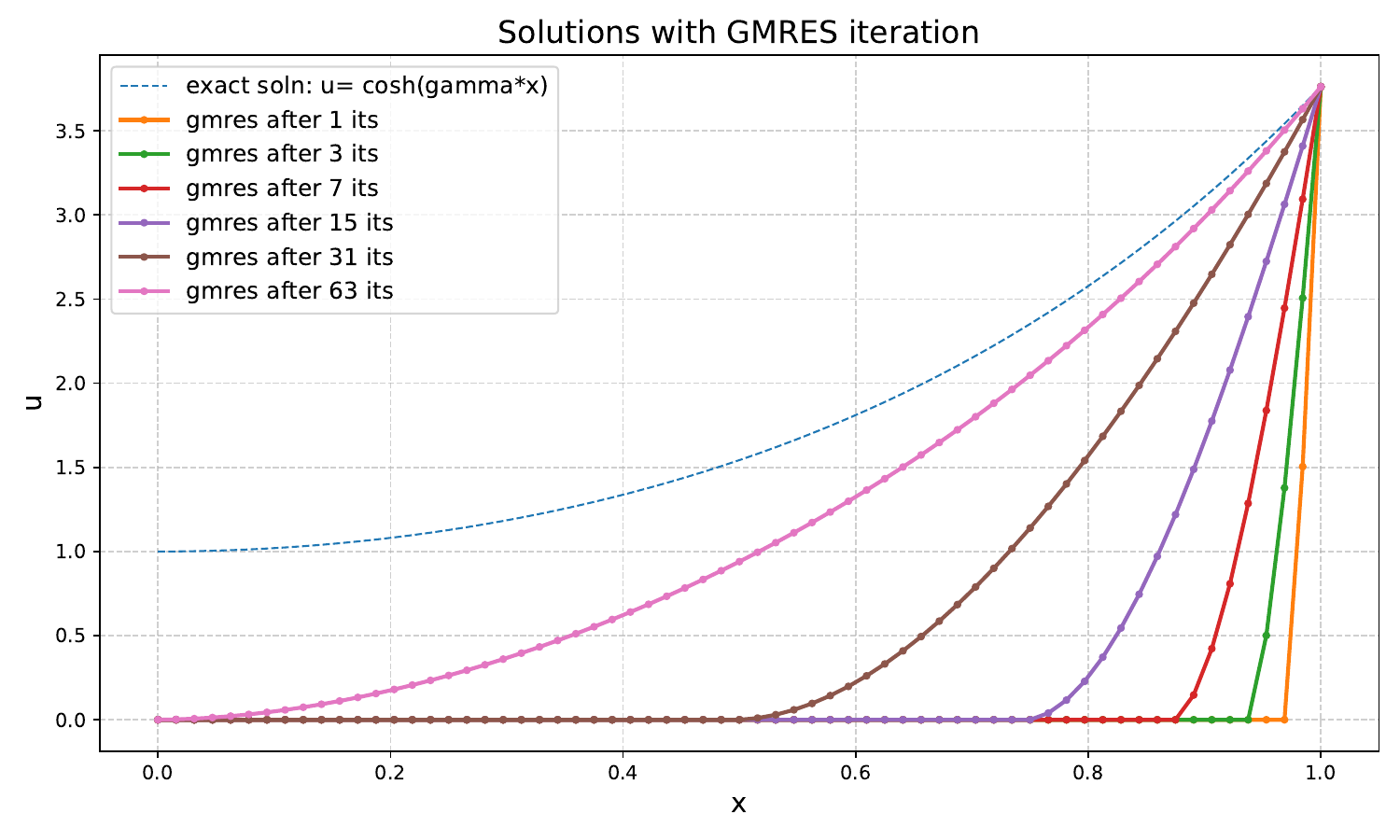}
    \caption{The true solution for eq.~(\ref{eq:example1}) compared with the results
    after $k= 1, 3, 7, 15, 31, 63$ \ac{GMRES} iterations.}
    \label{fig:gmres_stopping}
\end{figure}
In Fig.~\ref{fig:gmres_stopping} we display intermediate GMRES iterates in analogy to
Fig.~\ref{fig:stopping_64}.
While GMRES generates intermediate approximations that differ from those of the
\ac{CG} method, the qualitative behavior remains essentially the same.
In particular, the error at the left endpoint stagnates at zero up to the
$(n-1)$-st iteration.
Only when $k=n=64$ is reached does the method produce a qualitatively correct solution.
At this point, the discrete linear system is solved exactly up to roundoff errors,
while the discretization error, of course, persists.

A principal drawback of the \ac{GMRES} algorithm is the requirement to compute and store
an explicit basis of the Krylov space $\mathcal{K}_k(\Am,\bv)$.
The memory required to store the basis vectors of
may become the dominant limitation.
Consequently, \ac{GMRES} is commonly employed with restarts
\cite{saad1986gmres,baker2005technique}, where the Krylov space dimension is capped at a prescribed size.
Once this limit is reached, the method is restarted and a new Krylov space is constructed,
using the final iterate of the previous cycle as the initial guess.
It is well known that restarted \ac{GMRES} can exhibit subtle and sometimes counterintuitive
convergence behavior; see, for example, \cite{embree2003tortoise}.
In the present setting, we are particularly interested in how restarts influence convergence
and how they interact with the \loli.
\begin{remark}
For symmetric (possibly indefinite) matrices, the \ac{MINRES} method \cite{paige1975solution}
offers an alternative that exploits symmetry and relies on short three-term recurrences,
analogous to \ac{CG}, thereby avoiding the storage of a full Krylov basis.
We do not explore this further in this report.
\end{remark}

\begin{figure}[ht]
    \centering
    \includegraphics[width=0.7\textwidth]{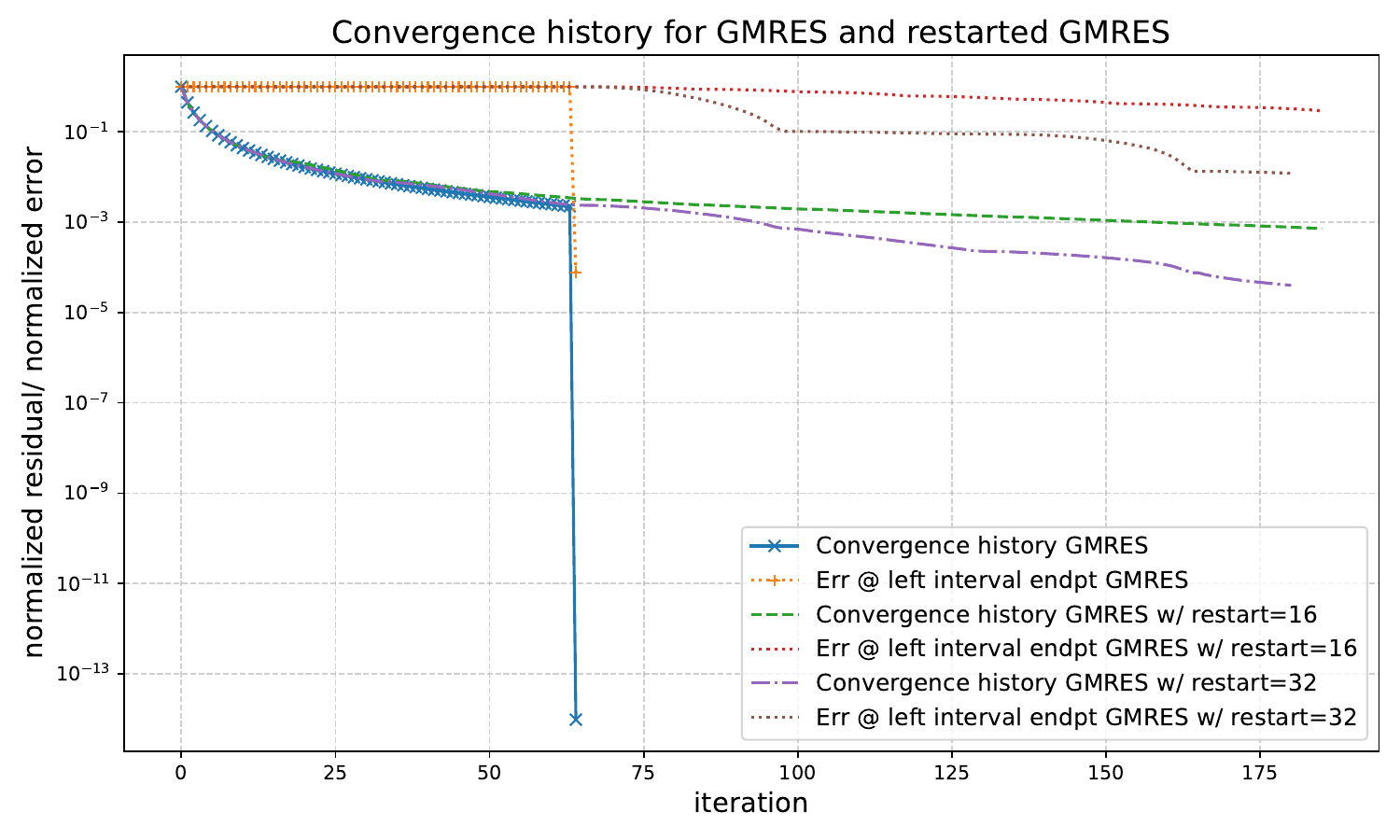}
    \caption{Convergence history for GMRES applied to eq.~(\ref{eq:example1}) with no restarts and restarts of 16 and 32, respectively.
    \label{fig:gmres_convergence}}
\end{figure}
As an initial investigation, we refer to Fig.~\ref{fig:gmres_convergence}, which reports the
convergence history in a format consistent with the previous figures.
The standard \ac{GMRES} algorithm without restarts exhibits the expected behavior:
the residual decreases slowly over most iterations and then drops abruptly by roughly ten orders
of magnitude in the final step, i.e., when $k=n$.
Similarly, the error at $x=0$ stagnates at the value $1$ throughout the iteration and
drops sharply at $k=n$ to approximately $10^{-4}$, corresponding to the discretization error.

Figure~\ref{fig:gmres_convergence} also reports the effect of employing restarts of size
$16$ and $32$.
In both cases, the characteristic final jump to the accurate solution is lost.
Instead, the convergence behavior becomes qualitatively similar to that observed for
stationary iterative schemes such as Gauss--Seidel or \ac{SOR}.
For $k>n$, both the residual and the pointwise error continue to decrease, but only slowly.
Within the iteration range shown ($k \leq 180$), neither restarted variant attains
discretization-error accuracy.
Evidently, restarting increases the total number of iterations required to reach
an accuracy commensurate with the underlying discretization.
For example, with restarts every $32$ iterations, the residual is reduced by roughly
four orders of magnitude after about $175$ iterations, whereas the error at $x=0$
decreases by only two orders of magnitude, which is still insufficient to reach
the discretization-error level.

\section{2D examples} \label{sec:advanced}
The primary application of the \ac{CG} method is the solution of large, sparse linear systems
that arise from the discretization of \acp{PDE} in two or three spatial dimensions.
A canonical example is the Poisson equation, which we adopt here as a model problem.
It is posed on the two--dimensional unit square $\Omega = (0,1)^2$ as
\begin{equation} \label{eq:poisson}
\begin{aligned}
 -\Delta u &= f \qquad &&\text{in } \Omega, \\
 u &= 0 \qquad &&\text{on } y=0 \text{ and } y=1, \\
 u &= \sin(\pi y)\cosh(\pi) \qquad &&\text{on } x=1, \\
 \frac{\partial u}{\partial x} &= 0 \qquad &&\text{on } x=0 .
\end{aligned}
\end{equation}
We emphasize that this problem is equipped with mixed Dirichlet and Neumann boundary conditions.
\begin{figure}[ht]
    \centering
    \includegraphics[width=0.75\textwidth]{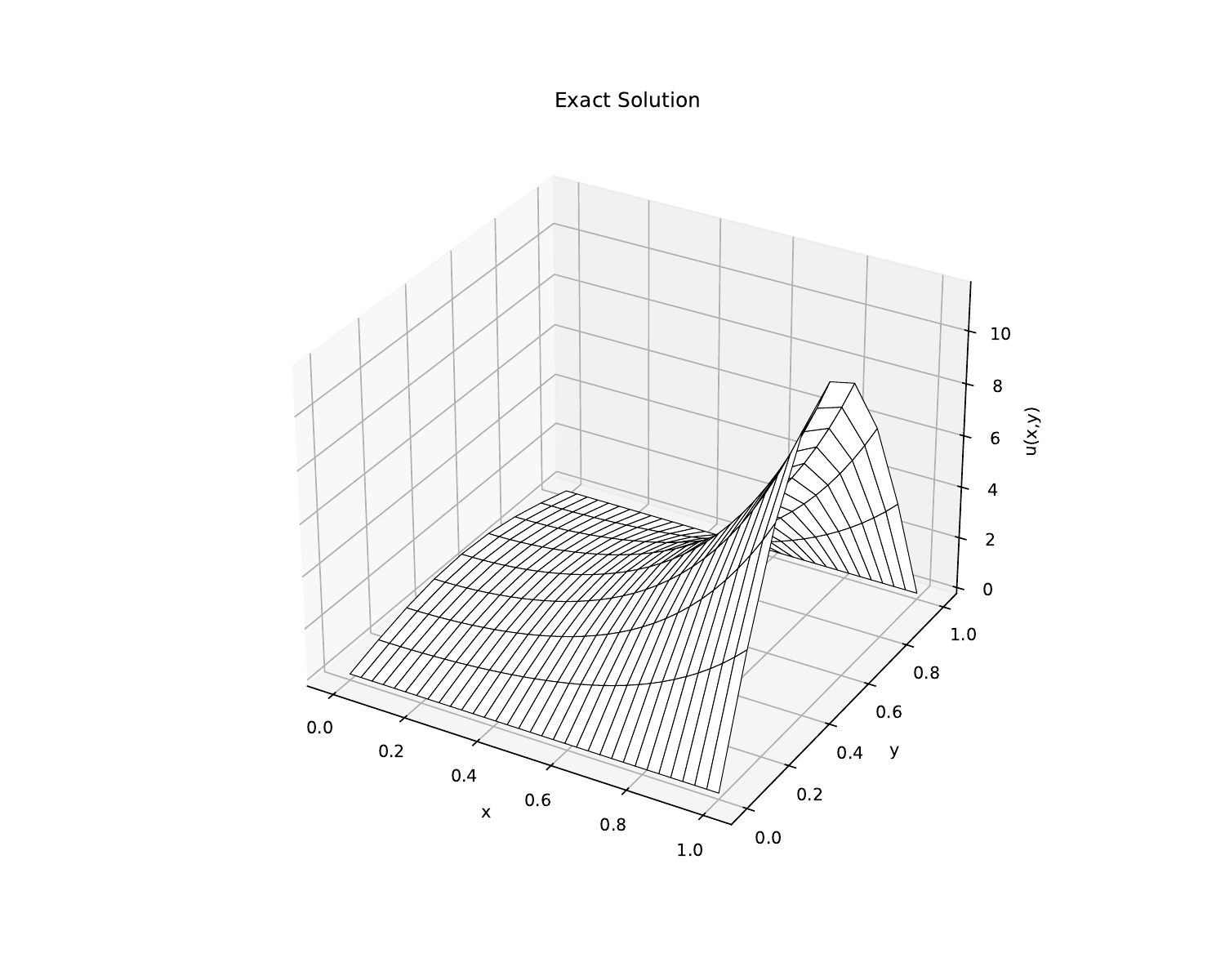}
    \vspace{-0.07\textwidth}
    \caption{Exact solution of  eq.~(\ref{eq:poisson}) displayed on a grid with $32 \times 8$ mesh nodes.
    \label{fig:poisson-exact}}
\end{figure} 
\begin{figure}[ht]
    \centering
    \vspace{-0.05\textwidth}	
    \includegraphics[width=0.75\textwidth]{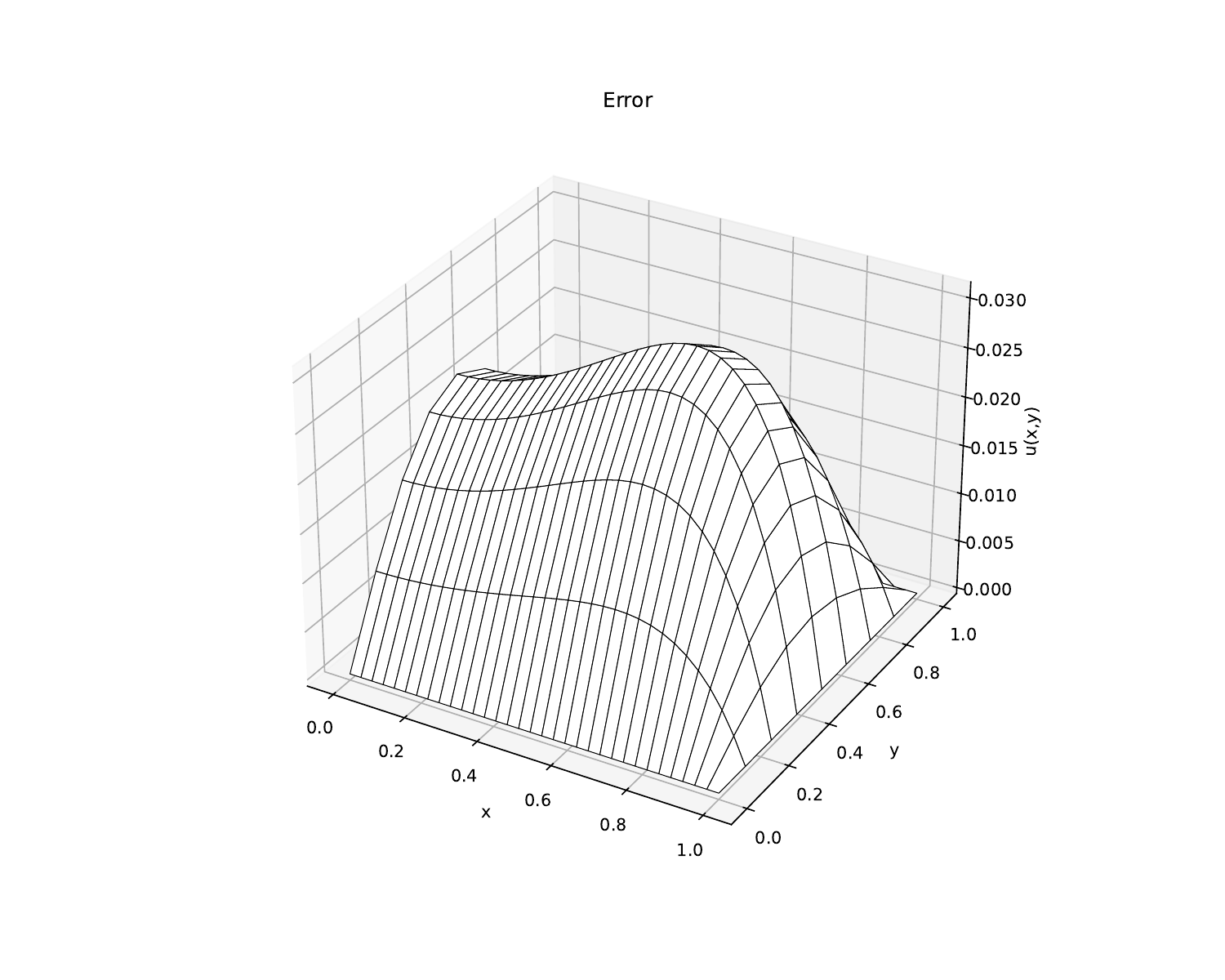}
    \vspace{-0.07\textwidth}
    \caption{Error of  discrete solution of eq.~(\ref{eq:poisson}) on a grid with $32 \times 8$ wrt.~to the exact solution.}
    \label{fig:poisson-error}
\end{figure} 
The attentive reader will notice the close analogy to~(\ref{eq:example1}).
Indeed, an analytical solution of~(\ref{eq:poisson}) obtained via separation of variables
leads to a one--dimensional subproblem of the form~(\ref{eq:example1}).
As a result, the exact solution of~(\ref{eq:poisson}) can be written explicitly as
\begin{equation}
  u^*(x,y) = \sin(\pi y)\,\cosh(\pi x).
\end{equation}

We discretize problem~(\ref{eq:poisson}) using the standard five--point finite difference stencil \cite{hackbusch2017elliptic} 
on a uniform Cartesian grid with $m$ and $n$ grid points in the $x$-- and $y$--directions, respectively.
For the isotropic case 
	$m = n \to \infty$, 
the scheme is second--order accurate.
In addition, we deliberately consider anisotropic grids with $m \neq n$ in order to investigate
their influence on the convergence behavior of the \ac{CG} method.
\begin{remark}
Up to a scaling factor, the same system matrix arises as the stiffness matrix of a conforming
finite element discretization with piecewise linear basis functions on triangular meshes,
when each rectangular grid cell is subdivided into two triangles along a diagonal; see, e.g.,
\cite{braess2007finite, elman2014finite, hackbusch2017elliptic}.
\end{remark}

The Neumann boundary condition is discretized in a manner that is consistent with the
second--order approximation used in the interior and such that it maintains symmetry.
After a priori elimination of the Dirichlet degrees of freedom, the resulting linear system is
symmetric positive definite and characterized by a sparse pentadiagonal structure.
For illustration, Fig.~\ref{fig:poisson-exact} depicts the analytical solution sampled on a uniform mesh with $32 \times 8$ cells,
while Fig.~\ref{fig:poisson-error} shows the corresponding error obtained by subtracting the
exact solution from the numerical approximation.

\begin{figure}[ht]
    \centering
    \vspace{-0.05\textwidth}	
    \includegraphics[width=0.75\textwidth]{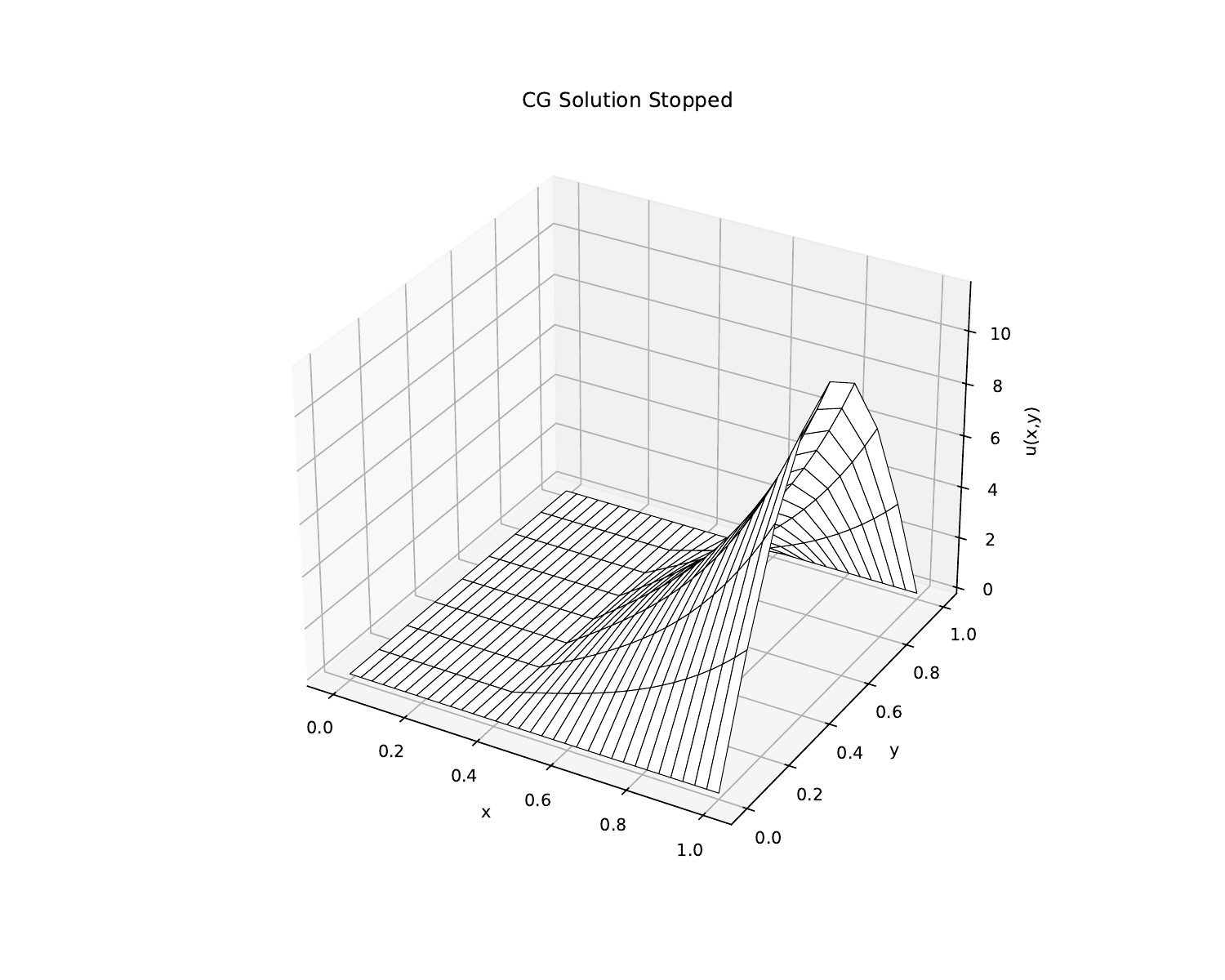}
    \vspace{-0.07\textwidth}
    \caption{Solution of discretized eq.~(\ref{eq:poisson}) on a grid with $32 \times 8$ mesh nodes with the
    	\ac{CG} method stopped after 19 iterations.}
    \label{fig:poisson-cg-stopped}
\end{figure} 
\begin{figure}[ht]
    \centering
    \includegraphics[width=0.7\textwidth]{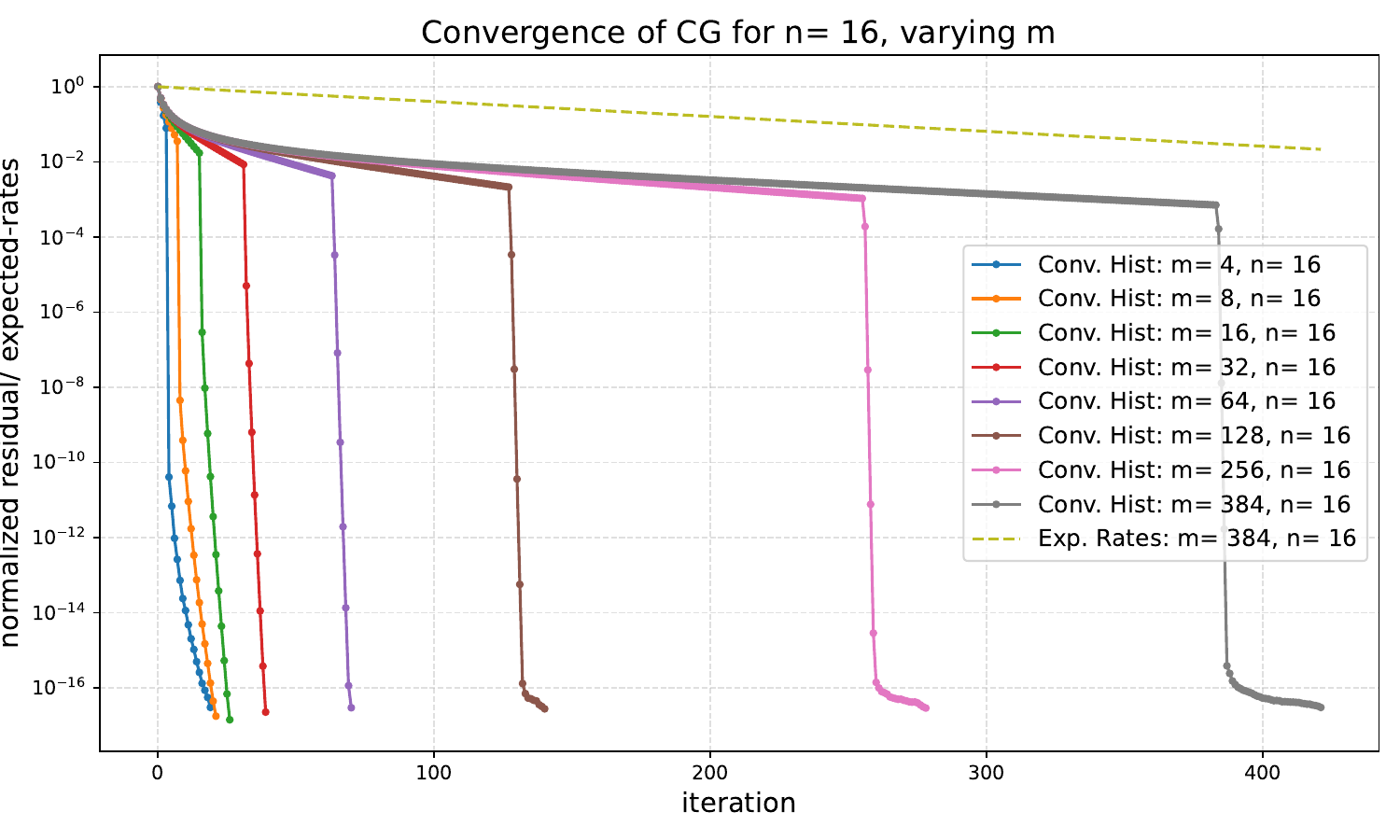}
    \caption{Convergence history for \ac{CG} applied to eq.~(\ref{eq:poisson}) on a grid with $m \times 16$ mesh cells, 
    	$m= 4, 8, 16, 32, 128, 256, 384$.}
    \label{fig:poisson-convergence-history}
\end{figure} 
We now apply the \ac{CG} method to solve the discretized form of~(\ref{eq:poisson}). 
Stopping the iteration after $19$ steps yields the approximation shown in Fig.~\ref{fig:poisson-cg-stopped}.
This figure clearly illustrates behavior analogous to the 1D case, highlighting the impact of the \loli:
in $19$ iterations, information from the Dirichlet boundary can propagate at most $19$ mesh cells.
Regions beyond this range remain essentially unaffected, preventing the \ac{CG} algorithm from producing a qualitatively meaningful solution.

A more detailed picture emerges from the convergence history of the residual, displayed in Fig.~\ref{fig:poisson-convergence-history}.
Here, $n=16$ is fixed, while $m$ is varied as indicated.
For reference, the expected asymptotic convergence rate is computed from the condition number
$\kappa(\Am) \approx 48139$, corresponding to the largest problem considered ($m=384$, $n=16$).
The residual histories for different mesh sizes exhibit behavior very similar to the 1D examples: 
an initial phase of slow, stagnating reduction persists until $k=m$ iterations have been completed.
Superlinear convergence sets in only after the \loli is overcome after $\min( m, n)$ iterations
and when information from the boundary has propagated across the entire domain.
Once this stage is reached, only a few additional iterations are required to reduce the residual to roundoff error.
Thus, for 2D (and by extension 3D) problems, the \loli again imposes a fundamental limitation on the performance of the \ac{CG} method.

\begin{remark}
In essentially isotropic cases, i.e., when $m \approx n$ in 2D, the impact of the \loli is less pronounced than in 1D. 
Here, the condition number scales as $\kappa(\Am) = \mathcal{O}(mn)$,
and thus~(\ref{eq:number_of_iterations}) predicts that the required number of iterations behaves like $k = \mathcal{O}(\sqrt{mn})$. 
Consequently, the additional constraint imposed by the locality limitation, $k > \min\{m,n\}$, 
does not exceed the iteration count that is suggested by the condition number, 
see, e.g.~eqn.~(\ref{eq:number_of_iterations}).
However, as illustrated in Fig.~\ref{fig:poisson-convergence-history}, 
this constraint can become significant for anisotropic meshes, 
where the disparity between $m$ and $n$ amplifies the effect of the \loli.
 \end{remark}
\begin{remark}
In terms of our analysis with anisotropic meshes to discretize the unit square, 
a quite similar situation would arise when using uniform meshes on a long and thin channel domain.
Furthermore, also in domains that are characterized by a system of connected thin channels, 
we must expect that they will each impose a severe \loli. This could, e.g., be the situation 
 in a geometrically resolved porous structure.
 \end{remark}
 \begin{remark}
The Poisson equation commonly arises from the continuity equation enforcing mass conservation in the
Stokes and incompressible Navier–Stokes equations. 
In such settings, prematurely terminating the \ac{CG} iteration while it remains in the stagnation phase—before information has fully propagated across the computational domain—can lead to violations of conservation properties. 
Consequently, the preservation of mass conservation is not ensured solely by the use of conservative discretization schemes; it also critically depends on the careful design of iterative solvers and on appropriate stopping criteria and error control within the solution process.
 \end{remark}
 
\section{Breaking the \loli by hierarchical preconditioning}
\label{sec:precond}
In this section we will attempt to accelerate the convergence of the \ac{CG} method by preconditioning.
In preconditioning, the usual goal is to find symmetry preserving transformations
of the form $\hat{\Am} = \Tm^T \Am \Tm$ such that the condition number $\kappa(\hat{\Am}) < \kappa(\Am)$.
In this situation, the convergence bound (\ref{eq:conv_bound}) will imply a faster convergence for the transformed system.

\begin{remark}
Note that the preconditioned \ac{CG} iteration is typically implemented without forming $\hat{\Am}$ explicitly,
since the \ac{CG} iteration on $\hat \Am$ can be realized
equivalently by multiplying the residuals $\rv_k$ with $\Cm = \Tm^T \Tm$ ,
see lines \ref{CG-residual0} and \ref{CG-residual} in Alg.~\ref{alg:cg}.
The symmetric matrix $\Cm$ is then referred to as the preconditioner.
For more details on preconditioning and its implementation see any of the standard references, 
such as, e.g, \cite{axelsson2001finite, saad2003iterative, axelsson2003iteration}.
\end{remark}

In this work, we adopt an alternative perspective on preconditioning, motivated by the
\loli of \ac{CG} methods.
Our goal is to identify transformations $\Tm$ and $\Cm$ that act as operators to accelerate
the propagation of information through the system.
The starting point is the elementary observation that a linear system with coefficient matrix $\Am$
can be transformed into an equivalent system by forming suitable linear combinations of its rows,
or, more generally, by replacing individual rows with linear combinations of other rows.

Such transformations can be described by left multiplication of $\Am$ with a nonsingular matrix.
While this process typically destroys the strict tridiagonal structure of $\Am$—at least partially—it enables the
\loli to be overcome.
This benefit comes at the cost of a partial loss of sparsity, as the transformed system
generally has an increased number of nonzeros per row and, consequently, potentially higher computational effort.

To illustrate the basic idea, we return to the 1D setting of example (\ref{eq:example1}).
We first consider the special case $n=8$ with mesh size $h=1/8$.
This choice is sufficiently small to allow an explicit presentation of the matrices involved,
while still large enough to reveal the essential structure of the proposed transformations.
As a first step, we introduce a transformation of the form
\begin{equation}
\Tm_8^4 =
\begin{pmatrix}
1 & \tfrac12 & 0 & 0 & 0 & 0 & 0 & 0 \\
0 & 1 & 0 & 0 & 0 & 0 & 0 & 0 \\
0 & \tfrac12 & 1 & \tfrac12 & 0 & 0 & 0 & 0 \\
0 & 0 & 0 & 1 & 0 & 0 & 0 & 0 \\
0 & 0 & 0 & \tfrac12 & 1 & \tfrac12 & 0 & 0 \\
0 & 0 & 0 & 0 & 0 & 1 & 0 & 0 \\
0 & 0 & 0 & 0 & 0 & \tfrac12 & 1 & \tfrac12 \\
0 & 0 & 0 & 0 & 0 & 0 & 0 & 1
\end{pmatrix}.
\end{equation}
Left multiplication by $\Tm_8^4$ replaces the rows with even indices $0,2,4,$ and $6$ by linear combinations
of their three neighboring rows, using weights $\tfrac12$, $1$, and $\tfrac12$, respectively,
while rows with odd indices remain unchanged.
This transformation can be interpreted as a single level of a one--dimensional hierarchical basis method;
see~\cite{yserentant1986multi} and~\cite{bank1988hierarchical}.

For the special case, when $\gamma=0$ and thus $d=2$ in eq.~(\ref{eq:1D-system}), the transformation results in
\begin{equation}
\Tm_8^4 \;  \Am(\tfrac18) = 
\frac{1}{h^2} \begin{pmatrix}
\frac{1}{2} & 0 & -\tfrac{1}{2} & 0 & 0 & 0 & 0 & 0\\
-1 & 2 & -1 & 0 & 0 & 0 & 0 & 0\\
-\tfrac{1}{2} & 0 & 1 & 0 & -\tfrac{1}{2} & 0 & 0 & 0\\
0 & 0 & -1 & 2 & -1 & 0 & 0 & 0\\
0 & 0 & -\tfrac{1}{2} & 0 & 1 & 0 & -\tfrac{1}{2} & 0\\
0 & 0 & 0 & 0 & -1 & 2 & -1 & 0\\
0 & 0 & 0 & 0 & -\tfrac{1}{2} & 0 & 1 & 0\\
0 & 0 & 0 & 0 & 0 & 0 & -1 & 2
\end{pmatrix}.
\end{equation}%
We note that in this case the number of non-zeros in each row is still at most 3 and that symmetry is lost.

The coefficients associated with rows $0,2,4,6$ (together with the final equation) now define a reduced system
for the corresponding unknowns.
Up to scaling and a reduction in dimension, this system retains the same structural form as the original one,
a fact that becomes apparent after reordering the equations by grouping even indices before odd ones.
For the given linear system, the transformation can also be interpreted as a single step of cyclic reduction,
in which all odd-indexed unknowns are eliminated, leaving a smaller system involving only the even unknowns.
This exact decoupling is specific to the present case with $d=2$.
For $d \neq 2$, the transformation no longer yields a perfect elimination; nevertheless, it can still be viewed
as an approximate elimination. 
In the resulting matrix, the remaining subsystems are then still weakly coupled.

Alternatively, the transformation could be generalized by adjusting the weights so that a cyclic elimination
is obtained even for $d\neq 2$.
For our purposes, however, the key observation is that, regardless of the diagonal value $d$, 
the transformation fundamentally alters the local coupling structure of the system.
In the transformed matrix, every second unknown is now directly coupled to unknowns at a distance of two grid points.
This extended coupling suggests a mechanism for mitigating the \loli identified in the preceding sections.
We investigate this effect with the following numerical experiments.
\begin{figure}[th]
    \centering
    \includegraphics[width=0.7\textwidth]{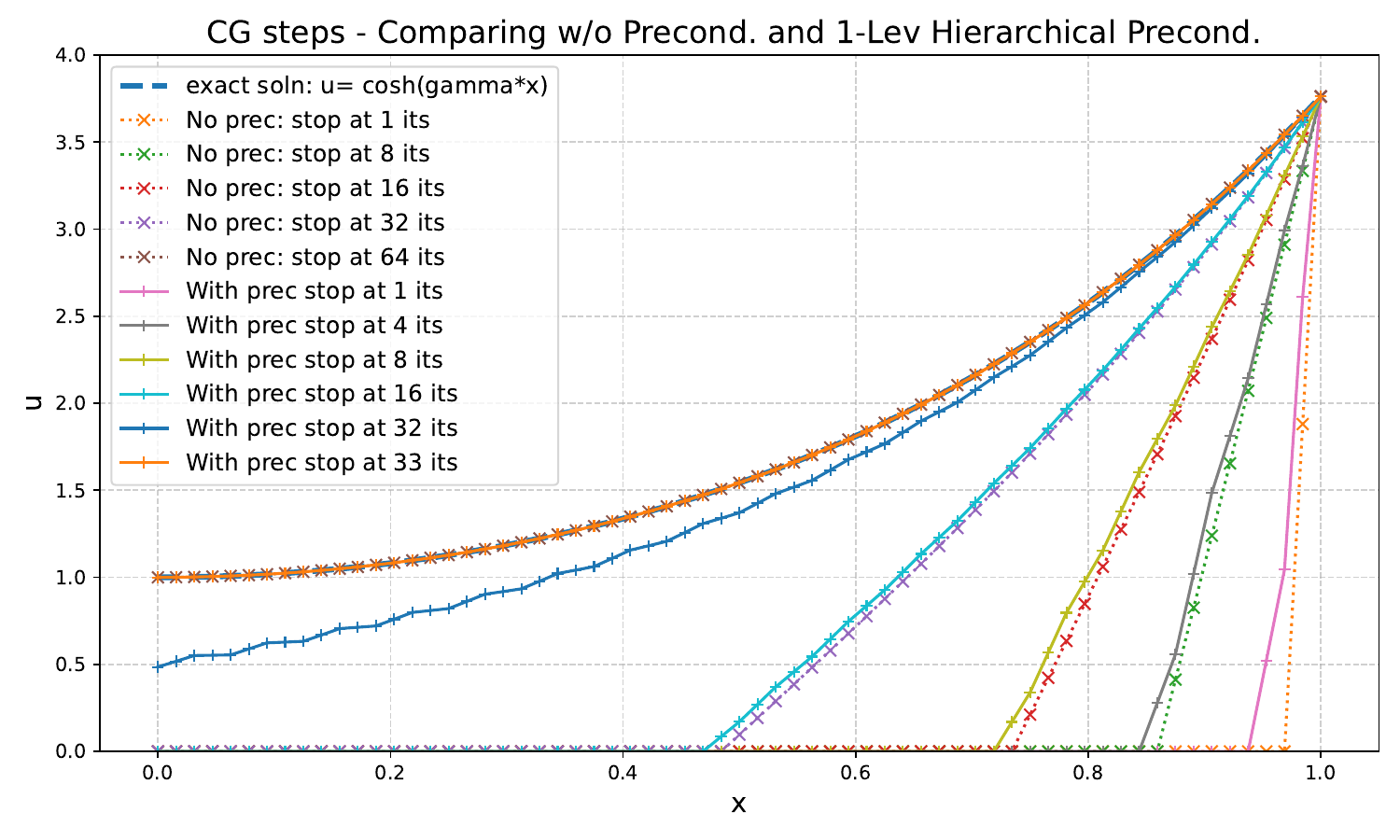}
    \caption{\ac{CG} with and without hierarchical 1-level preconditioner 
    for the solution of (\ref{eq:tridiag}) with $\gamma=2$ ($d$ accordingly).
    The iterates are displayed in comparison, showing that the preconditioned iteration delivers almost
    the same intermediate results in about half the number of iterations. A solution visually coinciding with the
    exact solution is reached in 33 iterations, rather than 64 for the unpreconditioned system.}
    \label{fig:precond-stopping-comparison}
\end{figure}

For use with the \ac{CG} iteration, we first symmetrize the transformation so that it becomes
\begin{equation}
\hat{\Am} := \Tm_8^4  \Am  (\Tm_8^4)^T = 
\frac{1}{h^2} 
\begin{pmatrix}
\tfrac{1}{2} & 0 & -\tfrac{1}{2} & 0 & 0 & 0 & 0 & 0 \\
0 & 2 & 0 & 0 & 0 & 0 & 0 & 0 \\
-\tfrac{1}{2} & 0 & 1 & 0 & -\tfrac{1}{2} & 0 & 0 & 0 \\
0 & 0 & 0 & 2 & 0 & 0 & 0 & 0 \\
0 & 0 & -\tfrac{1}{2} & 0 & 1 & 0 & -\tfrac{1}{2} & 0 \\
0 & 0 & 0 & 0 & 0 & 2 & 0 & 0 \\
0 & 0 & 0 & 0 & -\tfrac{1}{2} & 0 & 1 & 0 \\
0 & 0 & 0 & 0 & 0 & 0 & 0 & 2
\end{pmatrix}.
\end{equation}
Here, the transformed matrix $\hat{\Am}$ is perfectly decomposed into two independent subsystems
corresponding to the even and odd unknowns.
The odd, i.e., \emph{fine--grid}, unknowns give rise to a simple diagonal system,
whereas the even, \emph{coarse--grid}, unknowns form a subsystem that retains the same structural form as the original problem,
albeit with only half the number of degrees of freedom.
We note that this exact decoupling holds only for $d=2$; for $d\neq 2$ it is lost unless the transformation
is suitably modified.

To assess the impact of this construction as a preconditioner within the \ac{CG} iteration,
we return to the setting $n=64$ in~(\ref{eq:tridiag}) and again choose $\gamma=2$, so that $d>2$.
In this case, the preconditioner induced by $\Tm_{64}^{32}$ no longer yields an exact elimination
of the fine--grid unknowns, but rather an approximate one.
Applying the preconditioner
\begin{equation*}
	\Cm_{64}^{32} = (\Tm_{64}^{32})^{T} \,\Tm_{64}^{32}
\end{equation*}%
within the \ac{CG} algorithm produces the intermediate iterates shown in Fig.~\ref{fig:precond-stopping-comparison}.
A qualitatively satisfactory solution is now obtained after $33$ iterations, compared to $64$ iterations
in the unpreconditioned case.
Both the iteration counts and the intermediate solutions clearly demonstrate that the effective speed of information
propagation has been increased by approximately a factor of two.
As a result, the improved \loli leads to a convergence rate that is roughly doubled.

\begin{remark}
We further note that this speedup of the iteration by a factor of two cannot be explained by an improved condition number alone.
The condition number changes from $\kappa(\Am) \approx 2572$  to 
$\kappa(\Cm_{64}^{32}  \Am)  \approx 1295$. 
Thus, according to (\ref{eq:number_of_iterations}),
the number of iterations to reach the same error level should improve only by a factor of $\sqrt{2572 / 1295} \approx 1.4$.
Here, however, we observe a speedup of about $2$ rather than only $1.4$.
\end{remark}

\begin{figure}[ht]
    \centering
    \includegraphics[width=0.7\textwidth]{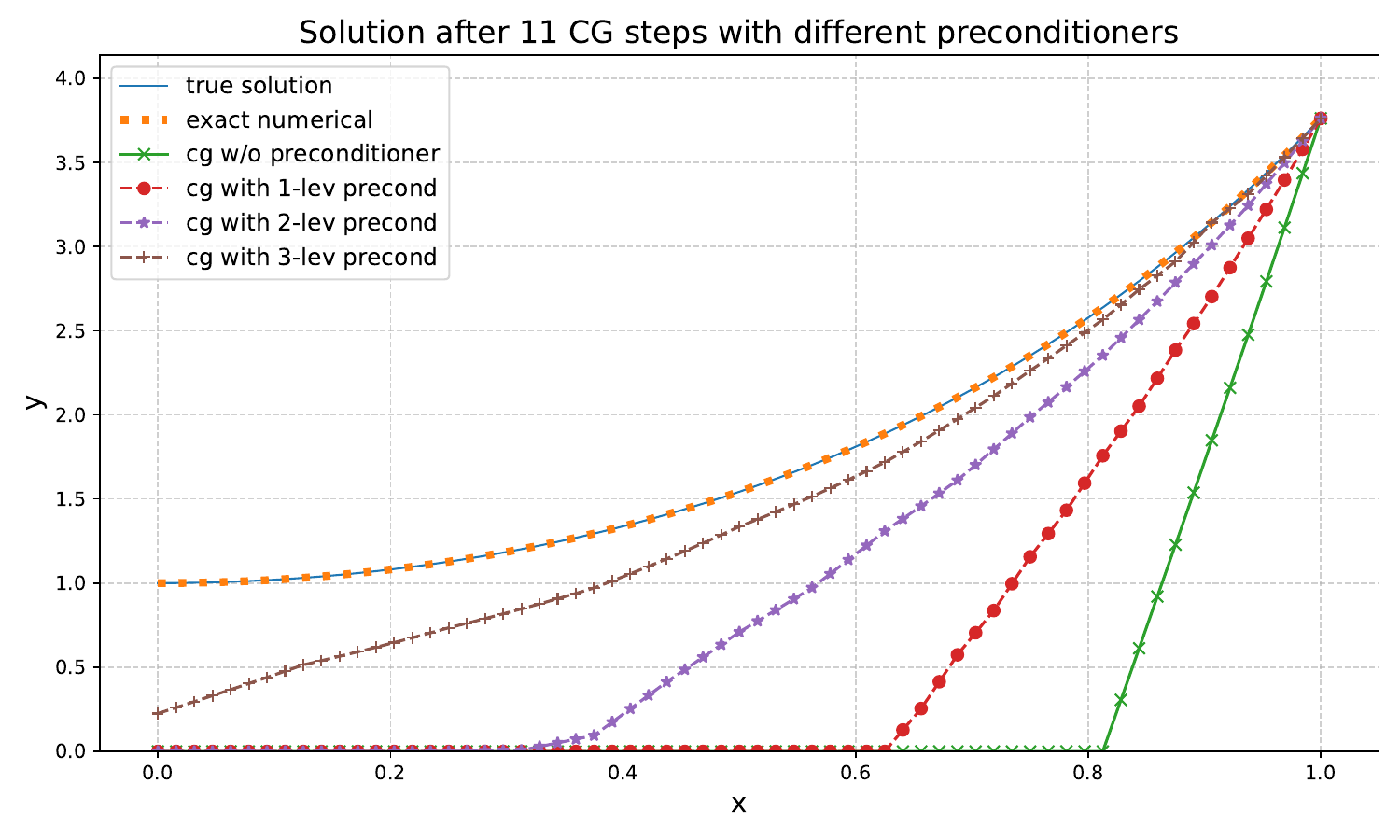}
    \caption{
    The iterates are shown after the 11th \ac{CG} iteration for different numbers of hierarchical preconditioning levels.
	Without preconditioning, the \loli restricts the propagation of information
	to 11 grid points away from the right boundary.
	Each additional level of preconditioning approximately doubles the effective propagation speed.
	With three levels of preconditioning, the approximation obtained after 11 iterations is still
	qualitatively inaccurate; however, only two further \ac{CG} iterations suffice to produce a solution
	that is visually indistinguishable from the exact one.
	\label{fig:precond-stopping-at-11}}

\end{figure}
Of course, it is now tempting to use the same idea recursively and apply it to the reduced system, 
so that a second, third and further levels of matrix rows  are identified and the correspondingly transformed. 
Clearly this will be most effective for larger systems, where this may then also be needed to achieve 
a satisfactory speed of convergence.

\begin{figure}[ht]
    \centering
    \includegraphics[width=0.7\textwidth]{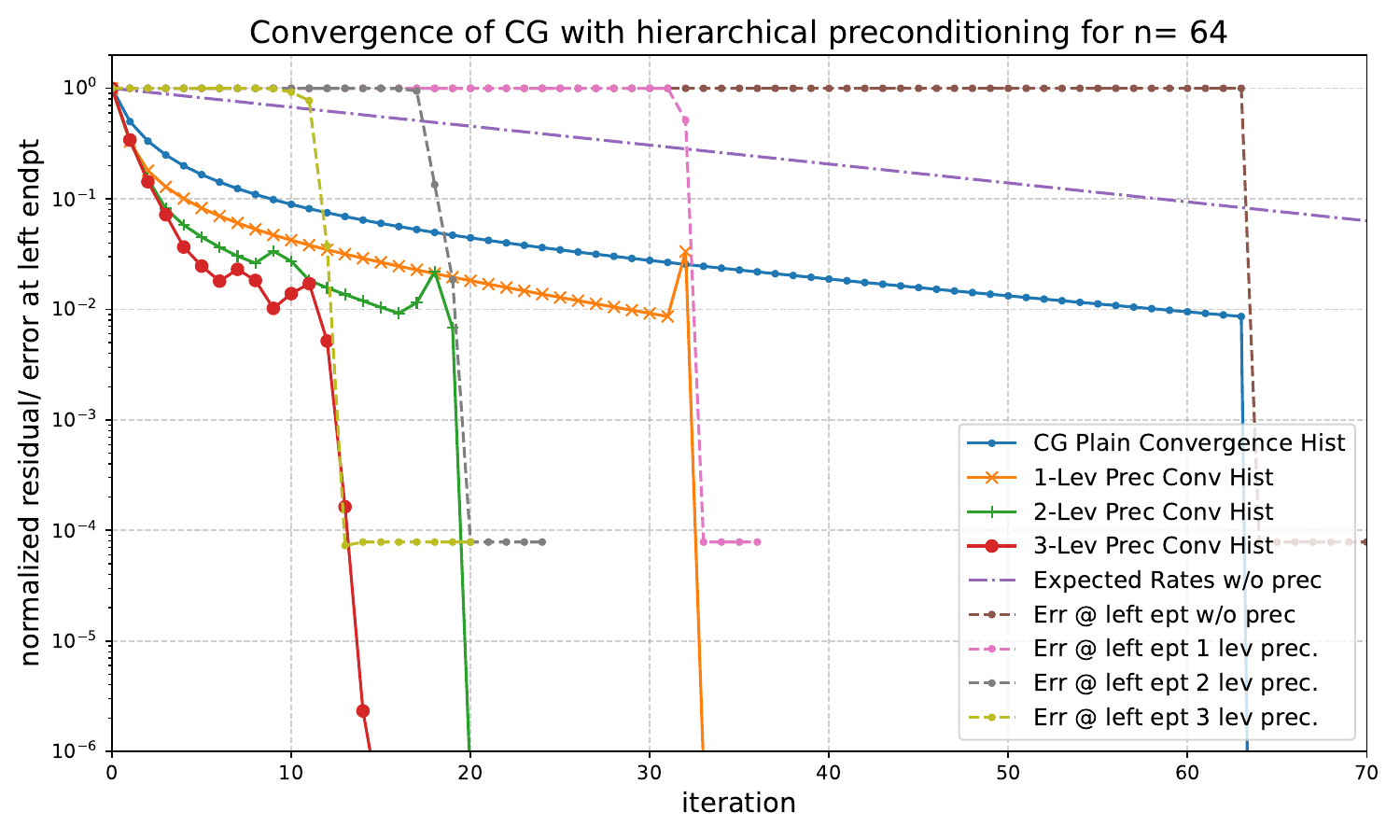}
    \caption{
    Convergence history for different levels of preconditioning and additionally the evolution of pointwise error at $x=0$.    
    \label{fig:prec-conv-hist}}
\end{figure}
We thus proceed by constructing such a preconditioner first using two hierarchical levels
and then also a further third level.
In Fig.~\ref{fig:precond-stopping-at-11} we display the solution after 11 iterations of the \ac{CG} method when these
different preconditioners are applied. 
With one or two levels of preconditioning, 
11 iterations are still not enough to propagate the information throughout the whole domain,
but the improvement becomes clearly visible.
With the 3-level hierarchical preconditioner, 11 iterations are sufficient so that the information
from the right boundary condition has reached also the leftmost point $x=0$ in the domain. 
However, the number of iterations is still not enough for the \ac{CG} method to find a qualitatively acceptable solution.
Nevertheless, two more iterations, 13 in total, are sufficient to achieve this.
Thus here the original iteration has not been accelerated by a factor of $2^3=8$ as might have been considered ideal.
However, with this, the acceleration is still by a factor of $64/13 \approx 5$.

\begin{remark}
The hierarchical preconditioning presented in this section is closely related to
hierarchical basis and multilevel preconditioners \cite{yserentant1986multi, oswald1994, ong1997hierarchical, ruede1993mathematical}
the BPX method \cite{bramble1990parallel}, the multi-level transformations of Griebel \cite{griebel1994multilevel, griebel1995parallel},
and various flavors of multigrid methods 
\cite{stueben1981multigrid, ruede1993fully, trottenberg2001multigrid, brandt2011multigrid, hackbusch2013multi, kuhn2022implicitly, bourne2023solver, litz2025memory, bohm2025large}.
\end{remark}

We next examine the convergence histories of the \ac{CG} method with the different preconditioners, 
both in terms of the residual norm and the pointwise error at the left endpoint of the interval ($x=0$).
In all cases, once the information has propagated across the entire domain
and the \loli is overcome, the residual exhibits rapid, superlinear convergence.
Each additional level of hierarchical preconditioning reduces the effective
\loli by approximately a factor of two.
The impact is even more pronounced when considering the pointwise error at $x=0$.
As in the unpreconditioned case, this error stagnates as long as the
\loli persists; however, once it is surpassed, the error decreases
rapidly and reaches the level of approximately $10^{-4}$, 
which corresponds to the remaining discretization error.

\section{Conclusions and Outlook}
\label{sec:conc}
In this technical report, we have presented an experimental investigation of the
\loli of iterative methods for the solution of sparse
symmetric positive definite linear systems.
Our experiments demonstrate that the efficiency of many iterative solvers is fundamentally constrained
by the slow propagation of information induced by the sparsity pattern of the system matrix.
In each iteration, information exchange is restricted to pairs of unknowns that are directly coupled
by nonzero matrix entries, a mechanism that is common to all Krylov subspace methods as well as to
classical stationary iterations.

Interpreting sparse matrices as graphs provides a natural conceptual framework for understanding this
phenomenon.
Each nonzero entry corresponds to an edge in the associated graph, and a single iteration propagates
information only to adjacent nodes.
As a consequence, the graph diameter yields a worst-case bound on the number of iterations required
before global information, such as boundary conditions, can influence the entire solution.
This observation explains the stagnation behavior frequently observed in early iterations and
highlights that the \loli represents a limitation that is independent of spectral considerations alone.

The numerical experiments presented in this work, primarily for one-dimensional model problems and
their extensions to two dimensions, illustrate that the \loli can dominate convergence behavior even
when the condition number indicates more favorable iteration counts.
In particular, mesh anisotropy exacerbates the effect by increasing the distance over which
information must be transported, thereby delaying the onset of rapid convergence.

These findings also provide an alternative motivation for multilevel methods and hierarchical
preconditioners.
From the perspective of the \loli, the principal role of multigrid techniques, hierarchical bases,
and related preconditioners is to accelerate the transport of information across the computational
domain by introducing couplings between distant unknowns.
Coarse-grid corrections and hierarchical transformations effectively reduce the graph diameter of
the preconditioned system, thereby relaxing the \loli and enabling faster global
information exchange.
Our experiments with hierarchical preconditioners clearly illustrate how successive levels
progressively alleviate the \loli and lead to substantial reductions in the number of iterations
required for convergence.

Beyond providing insight into the specific behavior of iterative solvers, 
this report also serves as an example of the systematic design of informative test problems and the extraction of mathematical understanding from computational experiments. 
This is an example of experimental mathematics, a scientific approach that remains comparatively
underdeveloped both in research practice and in educational settings.

Future work will focus on extending the present analysis to nonsymmetric and indefinite linear
systems, as well as on developing a graph-theoretic theory of the \loli.
Such a theory would provide a framework for understanding the effectiveness of multigrid,
domain decomposition, and preconditioning strategies from the standpoint of information propagation.

\section*{Acknowledgements}
The research was co-funded by the financial support of the European
Union under the REFRESH -- Research Excellence For Region Sustainability
and High-tech Industries project number CZ.10.03.01/00/22\_003/0000048
via the Operational Programme Just Transition.

The author gratefully acknowledges the scientific support and HPC resources provided by the Erlangen National High Performance Computing Center (NHR@FAU) of the Friedrich-AlexanderUniversit\"at Erlangen-N\"urnberg (FAU).
NHR funding is provided by federal and Bavarian state authorities.
NHR@FAU hardware is partially funded by the German Research Foundation (DFG), 440719683.

Portions of this manuscript were prepared and edited with the assistance of large language model–based tools. 
Similar tools were also used in the development of parts of the computational codes employed in this study. 
The author remains solely responsible for the scientific content, correctness, and interpretation of all results presented in this work.

\bibliographystyle{siamplain}
\bibliography{references}

\begin{thebibliography}{10}

\bibitem{arioli1992stopping}
{\sc M.~Arioli, I.~Duff, and D.~Ruiz}, {\em Stopping criteria for iterative
  solvers}, SIAM Journal on Matrix Analysis and Applications, 13 (1992),
  pp.~138--144.

\bibitem{Axelsson1994}
{\sc O.~Axelsson}, {\em Iterative Solution Methods}, Cambridge University
  Press, 1994.

\bibitem{axelsson2003iteration}
{\sc O.~Axelsson}, {\em Iteration number for the {Conjugate Gradient} 
  method}, Mathematics and Computers in Simulation, 61 (2003), pp.~421--435.

\bibitem{axelsson2001finite}
{\sc O.~Axelsson and V.~A. Barker}, {\em Finite Element Solution of Boundary
  Value Problems: Theory and Computation}, SIAM, 2001.

\bibitem{baker2005technique}
{\sc A.~H. Baker, E.~R. Jessup, and T.~Manteuffel}, {\em A technique for
  accelerating the convergence of restarted {GMRES}}, SIAM Journal on Matrix
  Analysis and Applications, 26 (2005), pp.~962--984.

\bibitem{bank1988hierarchical}
{\sc R.~E. Bank, T.~F. Dupont, and H.~Yserentant}, {\em The hierarchical basis
  multigrid method}, Numerische Mathematik, 52 (1988), pp.~427--458.

\bibitem{bohm2025large}
{\sc F.~B{\"o}hm, N.~Kohl, H.~K{\"o}stler, and U.~R{\"u}de}, {\em Large-scale
  multigrid with adaptive galerkin coarsening}, arXiv preprint
  arXiv:2511.13109,  (2025).

\bibitem{boisvert1991algorithms}
{\sc R.~F. Boisvert}, {\em Algorithms for special tridiagonal systems}, SIAM
  Journal on Scientific and Statistical Computing, 12 (1991), pp.~423--442.

\bibitem{bondeli1994cyclic}
{\sc S.~Bondeli and W.~Gander}, {\em Cyclic reduction for special tridiagonal
  systems}, SIAM Journal on Matrix Analysis and Applications, 15 (1994),
  pp.~321--330.

\bibitem{bourne2023solver}
{\sc E.~Bourne, P.~Leleux, K.~Kormann, C.~Kruse, V.~Grandgirard, Y.~Güçlü,
  M.~J. Kühn, U.~Rüde, E.~Sonnendrücker, and E.~Zoni}, {\em Solver
  comparison for {Poisson}-like equations on {Tokamak} geometries}, Journal of
  Computational Physics, 488 (2023), p.~112249,
  \url{https://doi.org/https://doi.org/10.1016/j.jcp.2023.112249}.

\bibitem{braess2007finite}
{\sc D.~Braess}, {\em Finite Elements: Theory, Fast Solvers, and Applications
  in Solid Mechanics}, Cambridge University Press, Cambridge, 3rd~ed., 2007.

\bibitem{bramble1990parallel}
{\sc J.~H. Bramble, J.~E. Pasciak, and J.~Xu}, {\em Parallel multilevel
  preconditioners}, Mathematics of Computation, 55 (1990), pp.~1--22.

\bibitem{brandt2011multigrid}
{\sc A.~Brandt and O.~E. Livne}, {\em Multigrid techniques: 1984 guide with
  applications to fluid dynamics, revised edition}, SIAM, 2011.

\bibitem{conte2017elementary}
{\sc S.~D. Conte and C.~De~Boor}, {\em Elementary numerical analysis: an
  algorithmic approach}, SIAM, 2017.

\bibitem{davis2016survey}
{\sc T.~A. Davis, S.~Rajamanickam, and W.~M. Sid-Lakhdar}, {\em A survey of
  direct methods for sparse linear systems}, Acta Numerica, 25 (2016),
  pp.~383--566.

\bibitem{duff2017direct}
{\sc I.~S. Duff, A.~M. Erisman, and J.~K. Reid}, {\em Direct methods for sparse
  matrices}, Oxford University Press, 2017.

\bibitem{elman2014finite}
{\sc H.~C. Elman, D.~J. Silvester, and A.~J. Wathen}, {\em Finite elements and
  fast iterative solvers: with applications in incompressible fluid dynamics},
  Oxford University Press, 2014,
  \url{https://doi.org/10.1093/acprof:oso/9780199678792.001.0001}.

\bibitem{embree2003tortoise}
{\sc M.~Embree}, {\em The tortoise and the hare restart {GMRES}}, SIAM Review,
  45 (2003), pp.~259--266.

\bibitem{greenbaum1997iterative}
{\sc A.~Greenbaum}, {\em Iterative Methods for Solving Linear Systems}, SIAM,
  1997, \url{https://doi.org/10.1137/1.9781611970937}.

\bibitem{greenbaum2021convergence}
{\sc A.~Greenbaum, H.~Liu, and T.~Chen}, {\em On the convergence rate of
  variants of the {Conjugate Gradient} algorithm in finite precision
  arithmetic}, SIAM Journal on Scientific Computing, 43 (2021), pp.~S496--S515.

\bibitem{greenbaum1996any}
{\sc A.~Greenbaum, V.~Pt{\'a}k, and Z.~e.~k. Strako{\v{s}}}, {\em Any
  nonincreasing convergence curve is possible for {GMRES}}, Siam Journal on
  Matrix Analysis and Applications, 17 (1996), pp.~465--469.

\bibitem{griebel1994multilevel}
{\sc M.~Griebel}, {\em Multilevel algorithms considered as iterative methods on
  semidefinite systems}, SIAM Journal on Scientific Computing, 15 (1994),
  pp.~547--565.

\bibitem{griebel1995parallel}
{\sc M.~Griebel}, {\em Parallel domain-oriented multilevel methods}, SIAM
  Journal on Scientific Computing, 16 (1995), pp.~1105--1125.

\bibitem{hackbusch2013multi}
{\sc W.~Hackbusch}, {\em Multi-grid methods and applications}, vol.~4 of
  Springer Series in Computational Mathematics, Springer, 1985,
  \url{https://doi.org/10.1007/978-3-662-02427-0}.

\bibitem{hackbusch2017elliptic}
{\sc W.~Hackbusch}, {\em Elliptic differential equations: theory and numerical
  treatment}, Springer, 2017.

\bibitem{HestenesStiefel1952}
{\sc M.~R. Hestenes and E.~Stiefel}, {\em Methods of {Conjugate Gradients} for
  solving linear systems}, Journal of Research of the National Bureau of
  Standards, 49 (1952), pp.~409--436.

\bibitem{hu2024solving}
{\sc X.~Hu and J.~Lin}, {\em Solving graph {Laplacians} via multilevel
  sparsifiers}, SIAM Journal on Scientific Computing, 46 (2024),
  pp.~S378--S400.

\bibitem{kuhn2022implicitly}
{\sc M.~J. K{\"u}hn, C.~Kruse, and U.~R{\"u}de}, {\em Implicitly extrapolated
  geometric multigrid on disk-like domains for the gyrokinetic {Poisson}
  equation from fusion plasma applications}, Journal of Scientific Computing,
  91 (2022), p.~28.

\bibitem{leveque2007finite}
{\sc R.~J. LeVeque}, {\em Finite difference methods for ordinary and partial
  differential equations: steady-state and time-dependent problems}, SIAM,
  2007.

\bibitem{litz2025memory}
{\sc J.~Litz, P.~Leleux, C.~Kruse, J.~Gedicke, and M.~J. K{\"u}hn}, {\em
  Memory-and compute-optimized geometric multigrid gmgpolar for curvilinear
  coordinate representations--applications to fusion plasma}, Journal of
  Computational and Applied Mathematics,  (2025), p.~117308.

\bibitem{malek2014preconditioning}
{\sc J.~M{\'a}lek and Z.~Strako{\v{s}}}, {\em Preconditioning and the
  {Conjugate Gradient}  method in the context of solving {PDEs}}, SIAM, 2014.

\bibitem{meurant1997computation}
{\sc G.~Meurant}, {\em The computation of bounds for the norm of the error in
  the {Conjugate Gradient}  algorithm}, Numerical Algorithms, 16 (1997),
  pp.~77--87.

\bibitem{meurant2024error}
{\sc G.~Meurant and P.~Tich{\`y}}, {\em Error Norm Estimation in the {Conjugate
  Gradient}  Algorithm}, SIAM, 2024.

\bibitem{ong1997hierarchical}
{\sc M.~E.~G. Ong}, {\em Hierarchical basis preconditioners in three
  dimensions}, SIAM Journal on Scientific Computing, 18 (1997), pp.~479--498.

\bibitem{oswald1994}
{\sc P.~Oswald}, {\em Multilevel finite element approximation: Theory and
  applications}, Teubner Skripten zur Numerik, B. G. Teubner, Stuttgart, 1994,
  \url{https://doi.org/10.1007/978-3-322-91215-2}.

\bibitem{paige1975solution}
{\sc C.~C. Paige and M.~A. Saunders}, {\em Solution of sparse indefinite
  systems of linear equations}, SIAM journal on Numerical Analysis, 12 (1975),
  pp.~617--629.

\bibitem{ruede1993fully}
{\sc U.~R{\"u}de}, {\em Fully adaptive multigrid methods}, SIAM Journal on
  Numerical Analysis, 30 (1993), pp.~230--248.

\bibitem{ruede1993mathematical}
{\sc U.~R{\"u}de}, {\em Mathematical and Computational Techniques for
  Multilevel Adaptive Methods}, SIAM, 1993.

\bibitem{saad2003iterative}
{\sc Y.~Saad}, {\em {Iterative methods for sparse linear systems}}, SIAM, 2003,
  \url{https://doi.org/10.1137/1.9780898718003}.

\bibitem{saad1986gmres}
{\sc Y.~Saad and M.~H. Schultz}, {\em {GMRES}: A generalized minimal residual
  algorithm for solving nonsymmetric linear systems}, SIAM Journal on
  Scientific and Statistical Computing, 7 (1986), pp.~856--869.

\bibitem{SpielmanCGDiameter}
{\sc D.~A. Spielman}, {\em The conjugate gradient and diameter}, in Spectral
  and Algebraic Graph Theory, Lecture Notes, Yale University, 2025,
  \url{http://cs-www.cs.yale.edu/homes/spielman/sagt/sagt.pdf}.

\bibitem{stueben1981multigrid}
{\sc K.~St{\"u}ben and U.~Trottenberg}, {\em Multigrid methods: Fundamental
  algorithms, model problem analysis and applications}, in Multigrid Methods,
  W.~Hackbusch and U.~Trottenberg, eds., Berlin, Heidelberg, 1982, Springer
  Berlin Heidelberg, pp.~1--176, \url{https://doi.org/10.1007/BFb0069928}.

\bibitem{press2007numerical}
{\sc S.~A. Teukolsky, B.~P. Flannery, W.~Press, and W.~Vetterling}, {\em
  Numerical Recipes 3rd edition: The Art of Scientific Computing}, Cambridge
  University Press, 2007.

\bibitem{trottenberg2001multigrid}
{\sc U.~Trottenberg, C.~W. Oosterlee, and A.~Schuller}, {\em Multigrid
  Methods}, Academic Press, 2001.

\bibitem{yserentant1986multi}
{\sc H.~Yserentant}, {\em On the multi-level splitting of finite element
  spaces}, Numerische Mathematik, 49 (1986), pp.~379--412.

\end{thebibliography}

\end{document}